%% file: main.tex
\documentclass[smallextended]{svjour3}       
\smartqed  
\usepackage{graphicx}
\usepackage[utf8]{inputenc} 
\usepackage[T1]{fontenc}
\usepackage{hyperref}
\usepackage{amsmath, amssymb}
\usepackage{mathtools}
\usepackage{enumerate,enumitem}
\usepackage{bm}
\usepackage{subcaption}
\usepackage{todonotes}
\usepackage{booktabs}
\usepackage{natbib}

\usepackage{tikz}
\usetikzlibrary{decorations} 
\usepackage{pgfplots}
\pgfplotsset{compat=1.8}

\usepackage[lined,algonl,boxruled,noend]{algorithm2e}

\DeclareMathOperator{\argmax}{argmax}

\newcommand{\adconP}{ \mathcal{Z}_P^a }
\newcommand{\longP}{ \mathcal{P}_Z(P) }
\newcommand{\adVP}{ \mathcal{T}_P^a }
\DeclareMathOperator{\init}{init}
\newcommand{\Vinit}{ V_{\init} }
\DeclareMathOperator{\pos}{pos}
\newcommand{\Path}{ (p_1,\dots ,p_\ell) }
\DeclareMathOperator{\ext}{ext}
\newcommand{\extp}{ p_{\ext} }
\DeclareMathOperator{\start}{start}
\newcommand{\startp}{ p_{\start} }
\newcommand{\tentacleP}{ \mathcal{T}_P }
\newcommand{\Dupl}{N_D}
\newcommand{\TT}{``Thurn \& Taxis''}
\newcommand{\G}{\ensuremath{G_{T}}}
\newcommand{\Gext}{\ensuremath{\widehat{G}_T}}

\begin{document}

\title{On the online path extension problem}
\subtitle{Location and routing problems in board games}


\author{Konstantin Kraus \and Kathrin Klamroth \and Michael Stiglmayr}


\institute{K.~Kraus, K.~Klamroth, M.~Stiglmayr \at
              University of Wuppertal, School of Mathematics and Natural Sciences, IMACM \\
              Tel.: +49-202-4393487\\
              \email{stiglmayr@uni-wuppertal.de}  
}

\date{Received: date / Accepted: date}

\maketitle
%

\begin{abstract}
    We consider an online version of a longest path problem in an undirected and planar graph that is motivated by a location and routing problem occurring in the board game \TT. Path extensions have to be selected based on only partial knowledge on the order in which nodes become available in later iterations. Besides board games, online path extension problems have applications in disaster relief management when infrastructure has to be rebuilt after  natural disasters. For example, flooding may affect large parts of a road network, and parts of the network may become available only iteratively and decisions may have to be made without the possibility of planning ahead.
    
    We suggest and analyse selection criteria that identify promising nodes (locations) for path extensions. We introduce the concept of \emph{tentacles} of paths as an indicator  for the future extendability. Different initialization and extension heuristics are suggested on compared to an ideal solution that is obtained by an integer linear programming formulation assuming complete knowledge, i.e., assuming that the complete sequence in which nodes become available is known beforehand. All algorithms are tested and evaluated on the original \TT\ graph, and on an extended version of the \TT\ graph, with different parameter settings. The numerical results confirm that the number of tentacles is a useful criterion when selecting path extensions, leading to near-optimal paths at relatively low computational costs.
    
\end{abstract}

\section{Introduction}\label{sec:intro}
\input{introduction}

\section{Brief review of related literature}\label{sec:lit}
\input{Literature}

\section{Online path extension problems}\label{sec:prob_def}
\input{ProblemAndDefinition}

\section{Solution methods}\label{sec:sol}
\input{SolutionMethods}

\section{Numerical results}\label{sec:res}
\input{Results}

\section{Conclusions and outlook}\label{sec:conclusions}

\input{ConclusionAndOutlook}

\section{Acknowledgements}
This work was partially supported by the project KoLBi (BMBF, Project-ID 01JA1507).


\end{document}

%% file: introduction.tex

\begin{figure}[htb]
    \small
    \centering
    \begin{tikzpicture}[scale=0.6,every node/.style={draw=black,circle,font=\footnotesize,inner sep=1pt,minimum size=15pt}]
    \node (1) at (0,0) {1};
    \node (2) at (2,0) {2};
    \node (3) at (4,0) {3};
    \node (4) at (6,0) {4};
    \node (5) at (8,0) {5};
    \node (6) at (0,-2) {6};
    \node (7) at (2,-2) {7};
    \node (8) at (4,-2) {8};
    \node (9) at (6,-2) {9};
    \node (10) at (8,-2) {10};
    \node (11) at (-2,-4) {11};
    \node (12) at (0,-4) {12};
    \node (13) at (2,-4) {13};
    \node (14) at (4,-4) {14};
    \node (15) at (6,-4) {15};
    \node (16) at (8,-4) {16};
    \node (17) at (10,-4) {17};
    \node (18) at (-2,-6) {18};
    \node (19) at (0,-6) {19};
    \node (20) at (2,-6) {20};
    \node (21) at (4,-6) {21};
    \node (22) at (6,-6) {22};
    \draw[-] (1) to (2);
    \draw[-] (1) to (6);
    \draw[-] (1) to (7);
    \draw[-] (2) to (3);
    \draw[-] (2) to (7);
    \draw[-] (3) to (4);
    \draw[-] (3) to (7);
    \draw[-] (3) to (8);
    \draw[-] (3) to (9);
    \draw[-] (4) to (5);
    \draw[-] (4) to (9);
    \draw[-] (4) to (10);
    \draw[-] (6) to (7);
    \draw[-] (6) to (11);
    \draw[-] (7) to (8);
    \draw[-] (7) to (12);
    \draw[-] (7) to (13);
    \draw[-] (8) to (9);
    \draw[-] (8) to (13);
    \draw[-] (8) to (14);
    \draw[-] (8) to (15);
    \draw[-] (9) to (15);
    \draw[-] (9) to (16);
    \draw[-] (10) to (17);
    \draw[-] (11) to (12);
    \draw[-] (11) to (18);
    \draw[-] (11) to (19);
    \draw[-] (12) to (13);
    \draw[-] (12) to (19);
    \draw[-] (12) to (20);
    \draw[-] (13) to (14);
    \draw[-] (13) to (20);
    \draw[-] (14) to (15);
    \draw[-] (14) to (20);
    \draw[-] (14) to (21);
    \draw[-] (15) to (16);
    \draw[-] (15) to (21);
    \draw[-] (15) to (22);
    \draw[-] (16) to (17);
    \draw[-] (16) to (22);
    \draw[-] (17) to (22);
    \draw[-] (18) to (19);
    \draw[-] (19) to (20);
    \draw[-] (20) to (21);
    \draw[-] (21) to (22);
    \end{tikzpicture}
    \caption{Graph $\G=(V_T,E_T)$ of the board of the game \TT. This graph has $n=|V_T|=22$ nodes and $m=|E_T|=45$ edges. \label{Fig:Boardgame as graph}}
\end{figure}
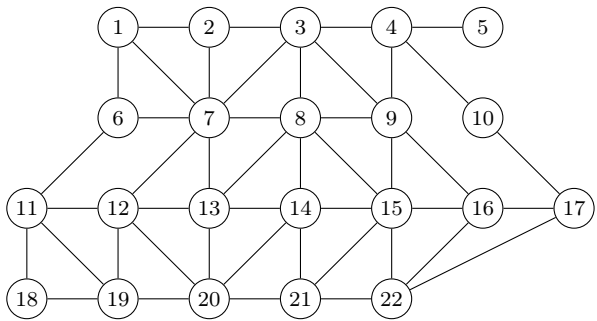

    Many board games are based on maps or networks and hence demand for locational decisions and/or routing solutions. Well-known examples are the games ``The Settlers of Catan'', where strategic locations for settlements are sought, and ``Ticket to Ride'' which asks for the construction of railroad networks under certain constraints. The \emph{online path extension problem} considered in this paper is motivated by the board game \TT, in which players aim to build the most rewarding postal carrier routes across Bavaria and surrounding regions. The board displays a map of selected cities (nodes) with a connecting road network (edges). It can be represented by the simple and planar graph $\G=(V_T,E_T)$ shown in Figure~\ref{Fig:Boardgame as graph}. Players build simple paths that represent postal carrier or taxi routes, one at a time, by consecutively collecting city cards that extend the currently active path. While aiming at longest paths in general, there is a trade-off between hoping for suitable cities to become available, and scoring (and hence closing) an intermediate path to reduce the risk of not being able to extend it any further, in which case the path is discarded. In this paper, we focus on the location and routing decisions that have to be made in this game, and hence ignore the game theoretic aspects of the game and omit further details regarding the rules. 

    The problem of selecting the next node from a finite candidate set of available nodes is modeled as a network location problem. Since the reward increases super-linearly with the length of a path, the driving criterion for node selection is the probability of being able to further extend a path in the subsequent iterations. This probability depends, on one hand, on the degree of the available nodes, and on the other hand on the number and type of nodes (i.e., city cards) that remain available in the deck, i.e., that will become available in future iterations. We introduce the concept of \emph{tentacles} of paths to predict the probability of long path extensions. In this context, a tentacle of a path is a node that is adjacent to (at least) one of the end nodes of the path while not already being a node in the path. We discuss modeling aspects as well as different solution approaches.
    
    The contribution of this paper is three-fold: First, we introduce a new and challenging class of online path extension problems. While being motivated by an application in the board game \TT, the problem has further applications of practical relevance. An example is the urgent repair of damaged infrastructure, including streets (edges) and shelter facilities (nodes), after natural disasters like flooding or earthquakes. Longer paths provide a wider operating distance, and path extensions are often only possible to nodes that are in some sense accessible. The information on accessible nodes may become available only iteratively, so that online path extensions have to be implemented without the possibility of planning ahead.
    
    Second, we develop heuristic path extension methods that are based on the novel concept of tentacles. Tentacles correspond to nodes that can be reached from the two end nodes of a given path and that are thus candidates for immediate path extensions. A path with a large number of (so-far) unused nodes that are adjacent to its end nodes has a high potential for further extensions. Thus, an extension towards a node that has many tentacles itself is preferable over an extension to a leaf node that allows for no future extensions.
    
    Third, we suggest an integer linear programming (IP) formulation for the online path extension problem under complete knowledge, i.e., assuming that the order in which further nodes become available is known beforehand. Paths constructed with this IP formulation are called \emph{ideal paths} and are used as a benchmark for the evaluation of the online path extension heuristics.

The paper is organized as follows: The literature on related problems is reviewed in Section~\ref{sec:lit}. A general problem formulation is given in Section~\ref{sec:prob_def}. Moreover, a formal definition of tentacles is provided. This section also contains an interpretation in the context of the board game \TT. An IP formulation for the computation of the ideal path under complete knowledge is provided in Section~\ref{sec:sol}. Moreover, four initialization heuristics and five online path extension heuristics are suggested in this section, all of which only work with the underlying graph and the partial information on available nodes in a given iteration. All algorithms are evaluated and compared in Section~\ref{sec:res} at randomly generated test instances (i.e., node orders) on the \TT\ graph shown in Figure~\ref{Fig:Boardgame as graph}, and on an extended version of the \TT\ graph. The paper is concluded with an outlook to further research directions in Section~\ref{sec:conclusions}.

%% file: Literature.tex
The literature on path extension problems in general and on selection criteria for promising nodes in the context of online path extension algorithms is scarce. 

Tentacles are related to leafs, and hence paths with many tentacles are in some sense related to maximum leaf spanning trees and the \emph{maximum leaf spanning tree problem} (MLSTP).
\cite{Fujie} used an integer programming approach to tackle the (MLSTP). Two different integer linear programming formulations are presented: An edge-vertex formulation and a vertex formulation. Moreover, valid inequalities are derived to strengthen the formulations. 
The exact algorithm of \cite{Fernau} operates on undirected graphs and is based on an equivalence to a related connected dominating set problem. 
\cite{Lu} suggest approximation algorithms for the MLSTP that start with an arbitrary tree. This tree is then improved by using $k$-changes (that switch $k$ tree edges with $k$ non-tree edges) until no further improvement is possible, yielding so-called $k$-locally optimal trees. They show that $k$-locally optimal trees approximate globally optimal trees with a high quality. They also present two heuristic algorithms for the MLSTP. 
\cite{Kneis} address the question whether a given directed graph contains an out-tree with at least $k$ leaves (where $k$ is considered an input parameter). They suggest a method that recursively grows a tree from a root node. This approach bears some similarity to the path extension procedure considered in this paper, however, aiming at directed graphs. Their algorithm improved the running time as compared to previous approaches. 
\cite{Reis} suggest a flow-based mixed-integer linear programming formulation for the MLSTP. They present numerical results that confirm the competitiveness of their model as compared to previous algorithms. We note that this flow-based formulation can be used as a basis for the identification of paths with fixed length and a maximum number of tentacles, which occurs as a subproblem in the context of this paper. 

The identification of paths with a large number of tentacles can also be related to so-called \emph{price collecting Steiner tree} (or path) problems. We exemplarily refer to \cite{Archer} who present improved approximation algorithms for the \emph{price collecting Steiner tree problem} (PCSTP) and also for the price collecting traveling salesman problem (PCTSP) and the price collecting path problem (PCPP). The goal is to find a tree (PCSTP), a cycle (PCTSP) or simple path (PCPP) that minimizes the overall costs \emph{and} the penalties for unused nodes.

Loosely related to the search for paths with many tentacles are approaches that aim at the generation of paths with specific properties. 
\cite{Khabbaz}, for example, study so-called \emph{heavy paths}. Their goal is to find a path of fixed length that maximizes the sum of the weights of the selected edges. They present an exact algorithm for the heavy path problem that uses a rank join approach. 
\cite{Awerbuch} investigates the multi-armed bandit problem where edge costs do vary over time. They design two randomized online adaptive routing algorithms for overlay networks. 

%% file: ProblemAndDefinition.tex
In this section, we introduce a mathematical terminology that facilitates the formulation of path extension problems. While some definitions and notation are tailored towards the board game \TT, most concepts are more generally applicable and are thus formulated in a general context.

\subsection{Paths, leafs, and tentacles}\label{subsec:tentac}


We first review some basic concepts from graph theory. For a general introduction into this topic, we refer to the textbooks \cite{ahuja93network,krumke12graphentheoretische}.
Let $G\coloneqq (V,E)$ be a simple and undirected graph with node set $V=\{v_1,\dots,v_n\}$, i.e., $|V|=n$, and edge set $E=\{e_1,\dots,e_m\}$, i.e., $|E|=m$. To simplify the notation, we will often refer to a node by its respective index, i.e., $v_i=i$, $i=1,\dots,n$.  Each edge $e_k\in E$ is defined by its two end nodes, i.e., $e_k=(v_i,v_j)=(i,j)$ with $v_i,v_j\in V$ and $v_i\neq v_j$.
Let $A=(a_{i,j})_{i=1,\dots ,n;\ j=1,\dots ,n}$ denote the adjacency matrix of $G$.
Then the \emph{degree} of a node $v_i\in V$ is given by $$\delta (v_{i})\coloneqq \sum_{j=1}^{n} a_{ij} .$$
We consider simple paths in $G$ that are defined as ordered sequences of pairwise different nodes. More precisely, $P=\Path $ is a \emph{simple path} of \emph{node length} $\ell(P)=\ell$ if $(p_i,p_{i+1})\in E$ for all $i\in\{ 1,\dots ,\ell-1\}$ and $p_{i}\neq p_{j}$ for all $i\neq j$. Note that we do not define the length of a path using the number of edges but the number of contained nodes. To emphasize this fact, we refer to it as the \emph{node length} of a path.

In this paper, we are interested in the iterative extension of simple paths by additional nodes with the ultimate goal of generating longest paths without node repetitions. In the context of the game \TT, such extensions have to be made based on partial knowledge on the set of available nodes since in each iteration only a subset of nodes is visible and additional nodes become available only in later iterations. Other possible applications of such \emph{online path extension problems} may occur, for example, when planning recovery operations after natural disasters that cut off the majority of roads in a transportation network, where roads (and thus nodes) become available only iteratively. Throughout this paper, we assume that (sub-)paths $P=(p_1,\dots,p_\ell)$ can only be extended at their end-nodes, i.e., at either of the two nodes $p_1$ or $p_\ell$. In such contexts, the \emph{prolongation potential} of a current (sub-)path plays a decisive role for the (iterative) path extension. Indeed, when selecting the next node such that, based on the partial knowledge available, \emph{many} further extensions of the path are known to be feasible in consecutive iterations, it is generally more likely that a long path can be constructed than when choosing a node that allows for only few extensions. This motivates the definition of \emph{tentacles} of a path.

\begin{definition}[Tentacles of a path]
    Let $G=(V,E)$ be a simple graph and let $P=\Path $ be a simple path in $G$. Then the \emph{tentacles} of $P$ are defined as $\tentacleP\coloneqq \{ v\in V\ :\ v\notin P \text{ and } ( (v,p_1)\in E \text{ or } (v,p_{\ell})\in E\}$. The \emph{number of tentacles} of $P$ will be denoted by \(\delta_P \coloneqq  |\tentacleP|\).
\end{definition}

\begin{figure}
    \small
    \centering
    \begin{tikzpicture}[scale=0.6,every node/.style={draw=black,circle,font=\footnotesize,inner sep=1pt,minimum size=10pt}]
    \node[dashed] (1) at (2,0) {};
    \node (2) at (4,0) {};
    \node (3) at (6,0) {};
    \node[dotted] (4) at (8,0) {};
    \node[dashed] (5) at (0,-2) {};
    \node (6) at (2,-2) {};
    \node[dashed] (7) at (4,-2) {};
    \node (8) at (6,-2) {};
    \node[dashed] (9) at (8,-2) {};
    \node[dotted] (10) at (4,-4) {};
    \node[dashed] (11) at (6,-4) {};
    \node[dotted] (12) at (8,-4) {};
    \node[dashed] (13) at (2,-4) {};
    \node[dotted] (14) at (4,-6) {};
    \node[dotted] (15) at (6,-6) {};
    \node[dotted] (16) at (4,2) {};
    \node[dotted] (17) at (6,2) {};
    \node[dotted] (18) at (8,2) {};
    \filldraw[black] (2) circle (10pt);
    \filldraw[black] (3) circle (10pt);
    \filldraw[black] (6) circle (10pt);
    \filldraw[black] (8) circle (10pt);
    \draw[-,thick] (2) to (3);
    \draw[-,thick] (2) to (6);
    \draw[-,thick] (3) to (8);
    \draw[dashed] (6) to (7);
    \draw[dotted] (7) to (10);
    \draw[dotted] (7) to (11);
    \draw[dotted] (10) to (11);
    \draw[dotted] (2) to (6);
    \draw[dashed] (8) to (11);
    \draw[dashed] (7) to (8);
    \draw[dotted] (1) to (2);
    \draw[dashed] (1) to (6);
    \draw[dotted] (2) to (16);
    \draw[dotted] (3) to (4);
    \draw[dotted] (3) to (17);
    \draw[dotted] (3) to (18);
    \draw[dashed] (5) to (6);
    \draw[dashed] (6) to (13);
    \draw[dashed] (8) to (9);
    \draw[dotted] (9) to (11);
    \draw[dotted] (10) to (13);
    \draw[dotted] (10) to (14);
    \draw[dotted] (11) to (12);
    \draw[dotted] (11) to (14);
    \draw[dotted] (11) to (15);
    \end{tikzpicture}
    \caption{A simple graph $G=(V,E)$, where $V$ comprises all dotted, dashed and filled nodes, with a simple path $P$ (filled nodes). The dashed nodes are the tentacles of $P$ and hence $\delta_P=6$.\label{Ex:TentacleDef}}
\end{figure}
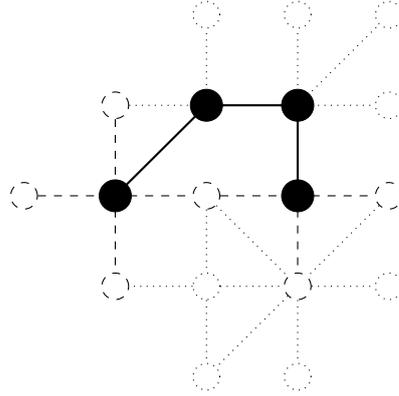

See Figure~\ref{Ex:TentacleDef} for an illustration of a simple path and its tentacles.

\begin{lemma}\label{lemma:bounds}
    Let $G=(V,E)$ be a simple graph and let $P=\Path$ be a simple path in $G$ with node length $\ell(P)=\ell$. Then the following properties hold:
    \begin{enumerate}
        \item $\delta_P\leq |V|-\ell$
        \item If $\ell=1$, then $\delta_P=\delta(p_1)$.
        \item If $\ell\geq 2$, then $\delta_P\leq \delta(p_1)+\delta(p_{\ell})-2$. This bound is tight when $G$ is a tree.
        \item If $\ell\geq 2$, then $\delta_P \geq \max\{\delta(p_1),\delta(p_\ell)\}-\ell$.
    \end{enumerate}
\end{lemma}
\begin{proof}
    The statements follow directly from the definition.
\end{proof}

Note that the bound in Lemma~\ref{lemma:bounds} (4) is also tight, for example, when $\delta(p_1)=\delta(p_\ell)$, both nodes $p_1$ and $p_\ell$ are adjacent to all nodes in $P$, and all other nodes that are adjacent to $p_1$ are also adjacent to $p_\ell$ (and vice versa).

\subsection{Available nodes and extension tuples}

In this section we focus on the particular setting in the game \TT, which may be interpreted as an online path extension problem. Here, path extensions can only be chosen from a finite selection of available nodes that has constant size during the course of the game. Such selections will be stored using vectors, or tuples, of nodes.  Selections may contain duplicates, i.e., nodes may occur several times in a selection, which reduces the number of choices. Whenever a node is selected as an extension, another node is added to the selection and hence becomes  available for the next iteration.
The three following definitions are required to mathematically describe one instance of the resulting game.

\begin{definition}[Initial tuple and initial set]
    Let $\Dupl\in\mathbb{N}$ denote the \emph{number of duplicates of every node}. Moreover, let $(i_1,\dots,i_c)$ with $i_k\in V$, $k=1,\dots ,c$, be an \emph{initial tuple} of \emph{size} $c\in\mathbb{N}$, referred to as \emph{initial $c$-tuple} in the following. Note that $i_l = i_k,\ l\neq k$ is not forbidden. 
    $\Vinit\coloneqq \lbrace i_1,\dots ,i_c \rbrace$ is called \emph{initial $c$-set}, without counting for duplicates.
\end{definition}

The initial $c$-set specifies the set of nodes that are available in the first iteration of the game. Hence, the first node of the path has to be selected from the set $\Vinit$. Moreover, we assume that all nodes in $V$ occur exactly $\Dupl$ times in a deck of cards (representing nodes) that is 
fixed from the beginning of the game, however, unknown to the players. This node deck can be represented by a vector of length $n\cdot \Dupl$, that results from combining an initial $c$-tuple with an appropriate \emph{extension tuple}. In this context, the extension tuple contains the information on possible path extensions during later iterations of the game.
    
\begin{definition}[Extension tuple]
    Let $d\coloneqq n\, \Dupl-c$ denote the total number of remaining cards in a given node deck that are not contained in the initial $c$-tuple. Then  $(j_1,\dots ,j_d)$ with $j_k\in V$, $k= 1,\ldots, d$, is called an \emph{extension tuple} of size $d$, also referred to as \emph{ordered extension $d$-tuple}, if the total number of duplicates is bounded by $\Dupl$ for every node $i\in V$, $i=1,\dots,n$. In other words, 
    $j_k=j_s, \; k\neq s$ is not forbidden, but there can be at most $\Dupl$ copies of the same node in 
    \((i_1,\ldots,i_c,j_1,\ldots,j_d)\), i.e., a feasible combination of initial tuple and extension tuple.
\end{definition}

To avoid confusion, we will refer to nodes from an initial tuple by $i_k$, $k\in\{1,\dots,c\}$, while the nodes from an extension tuple will be denoted by $j_k$, $k\in\{1,\dots,d\}$.
Combining an initial $c$-tuple with an extension tuple thus yields a deck of exactly $n\cdot \Dupl$ cards. Note that, while the order of nodes in an extension tuple is decisive for the order in which new nodes become available during the course of the game, the order of nodes is not relevant for the initial $c$-tuple since all of its nodes become available simultaneously at the initialization of the game. An instance of the game can now be defined using \emph{feasible settings} that combine an initial $c$-tuple with an appropriate extension tuple such that each node or card, respectively, occurs exactly $\Dupl$ times.

\begin{definition}[Feasible setting]
    A \emph{feasible setting} consists of an initial $c$-tuple and an ordered extension $d$-tuple, such that $d=|V|\cdot \Dupl-c$ holds and such that there are exactly $\Dupl$ copies of every node in the combination of both tuples.
    A \emph{$(c,\Dupl)$-setting} is a family of all feasible settings for a given problem size. 
\end{definition}

Note that a feasible setting could alternatively be defined as a node deck, i.e., a vector with $n\cdot\Dupl$ components from the set $V$ such that each element of $V$ occurs exactly $\Dupl$ times. The distinction between the initial tuple and the extension tuple, however, will turn out useful in the context of the online path extension problem considered in this paper.

\subsection{Online path extension problem}\label{subsec:OPEP}

Now we can describe an instance of the online path extension problem. Assume that a simple graph and a feasible $(c,\Dupl)$-setting are given consisting of an initial tuple $(i_1,\dots,i_c)$ and and extension tuple $(j_1,\dots,j_d)$. Then we want to find a simple path $P$ of maximum node length that is generated iteratively according to the following \emph{online path extension procedure} (OPEP):
\begin{enumerate}
    \item Select one node $\startp\in\Vinit$ from the initial tuple. Set $k=1$ and $P_0=(\startp)$.
    \item Select one node $\extp\in\Vinit\cup\lbrace j_1,\dots,j_k \rbrace$ which is adjacent to at least one end node of $P_{k-1}$. If there is no such node: Stop.
    \item Add $\extp$ at the end or at the beginning of the path $P_{k-1}$ (such that it is adjacent to that endpoint) to create $P_k$. Increase $k$ by one.
    \item Continue with (2) or end the procedure.
\end{enumerate}

Note that during this procedure, the extension tuple is generally only partly known. Indeed, in step (1) there is no information about the order of the extension tuple and in step (2) only the information on the first $k$ elements (where $k$ is the node length of the current path) of the extension tuple is known.

\begin{definition}[TT-path]
    Let a path $P=\Path$ and a feasible $(c,\Dupl)$-setting with initial tuple $(i_1,\dots ,i_c)$ and extension tuple $(j_1,\dots ,j_d)$ be given. $P$ is a \emph{TT-path (\TT-path)} w.r.t.\ the given feasible setting if and only if it can be generated by the online path extension procedure.
\end{definition}

\subsection{Interpretation in terms of the game ``Turn and Taxis''}

One application of an online path extension problem is a simplified one-player version of the board game \TT. In this game, players operate on a graph with $22$~nodes which represent different cities in southern Germany, Austria, Switzerland, Czechia and Poland. The game comes to an end whenever one player manages to build $20$~postal offices, under the constraint that there can not be more than one postal office per player per city. In order to build postal offices, players iteratively extend paths (that represent potential postal routes). In each iteration, a player has to add an additional node to his or her current path. The player can then decide to ``close'' this path and build postal offices in all of its nodes, or to try to further extend the path in the next round. The latter choice bears some risk: If it turns out impossible to extend the path in the next round, the path has to be discarded and no postal offices can be build on any of its nodes. The player then has to start over with a new path. There are always six cards on display from which a player can select one node to start or to extend his or her path. To extend a path, the next node must be added to the first or to the last node of the current path.

The number of duplicates indicates how many cards of every city are in the set of cards. There are three copies of every city in the set of cards of the original game \TT, i.e., $\Dupl=3$ in the original setting. A feasible setting corresponds to a complete stack of cards, divided into a part that is visible at the beginning of a game (the initial tuple) and stack of cards 
that is drawn as the game progresses (the extension tuple). We aim to heuristically model the decision process of a single player variant of \TT.

In the first iteration of a single-player version of this game, the player has to choose a card of a city from the initial tuple to initialize a path while not having any information on the order of the cards in the extension tuple. In all later iterations, the player can select a card from the initial tuple or from the first $k$ elements of the extension tuple, where $k$ is the total node length of the path computed so far. Again, no information on the order of the remaining cards in the extension tuple (apart from the first $k$ cards) is available to the player.

In this paper, we focus on the generation of the \emph{first} path and stop the path construction whenever no additional node can be added to one of its end nodes. We consider a single-player version in the sense that there is no competition, i.e., there are no other players competing about the available cards. Moreover, we evaluate the final path node length without requiring the player to make a stopping decision.

%% file: SolutionMethods.tex
The main focus of our work are heuristic path initialization and path extension strategies that iteratively select nodes based on the partial information available in a given iteration, according to the online path extension procedure. The quality of these strategies is evaluated in comparison with benchmark solutions that are computed based on complete knowledge on the order of nodes in the considered feasible setting. Such a benchmark solution can be computed by using an appropriate integer programming (IP) formulation. The solution to this IP formulation is referred to as an \emph{ideal} solution in the following, since it assumes complete information on the initial tuple and the extension tuple (which is actually not available to the player). We will compare the length of heuristically computed paths to the length of the ideal solution to evaluate the quality of the different heuristics.

\subsection{Exact solution assuming complete knowledge -- an IP model}\label{subsec:exsol}
\input{IP}

\subsection{Path initialization and (online) path extension heuristics based on partial knowledge}\label{subsec:heu}
\input{heuristics}

%% file: IP.tex
Let a simple graph $G=(V,E)$, its adjacency matrix $A=(a_{ij})$, and a feasible $(c,N_D)$-setting with $d=n\cdot N_D-c$ be given, with initial set $\Vinit$ and an extension tuple $(j_1,\dots ,j_d)$.
Based on this data, we define an \emph{availability matrix} $D\in\lbrace 0,1 \rbrace^{(d+1) \times n }$ that contains the information on the order in which nodes become available. Thus, $D_{k,p}=1$ if and only if a card that represents the node $p$ is revealed in iteration $k$. More precisely, we set
\begin{align*}
    D_{1,p} & =
\begin{cases}1 & \text{if}\; p\in \Vinit \\ 0 &  \text{otherwise} \end{cases} && \quad \forall p\in V, \\
    D_{k,p} & =\begin{cases}1 &\text{if}\;  p=j_{k-1}  \\ 0 &  \text{otherwise} \end{cases} && \quad \forall p\in V,\ k\in\lbrace2,\dots,d+1\rbrace.
\end{align*}
Note that the first row of the availability matrix $D$ contains at most $c$ non-zero entries (since $c$ nodes are available in the initial $c$-tuple and duplicate nodes are not counted in the availability matrix). The other rows of $D$ contain exactly one non-zero entry since only one additional node becomes available per iteration. 

For a given instance defined by the corresponding availability matrix $D$, we now aim at finding a longest TT-path, i.e., a longest path $P$ that is constructed according to the online path extension procedure introduced in Section~\ref{subsec:OPEP}. In this section, we assume complete knowledge of the matrix $D$.
Towards this end, let
the index $k\in\{ 1,\dots ,d+1\}$ denote an iteration, and let the index $p \in V$ denote a node. We introduce 
two sets of binary decision variables: $x_{kp}\in\{0,1\}$ indicates whether node $p$ is chosen in iteration $k$ or not, 
and $E_{kp}^{s}\in\{0,1\}$ and $E_{kp}^t\in\{0,1\}$ indicate whether node $p$ is the start node or the end node of the path in iteration $k$, respectively. 
Since at most one node can be chosen in each iteration, we want to maximize the number of iterations in which an additional node can be feasibly chosen according to the online path extension procedure. Whenever a further extension of the current path is impossible in an iteration $\bar{k}$, all variables $x_{kp}$ are equal to zero for all consecutive iterations $o\geq\bar{k}$. As a consequence, the node length of the final path is given by the sum over all decision variables $x_{kp}$. This leads to the following integer linear programming formulation:

{\small\begin{subequations}
\begin{align}
    \max \quad & \sum_{k=1}^{d+1} \sum_{p=1}^{n} x_{kp} \label{eq:1}\\[2ex] 
    \text{s.t.} 
    \quad & x_{kp} \leq \sum_{\kappa=1}^{k} D_{\kappa p} & \forall k\in\{1,\dots ,d+1\},\ p\in V \label{eq:2}\\
    & \sum_{\kappa=1}^{d+1} x_{\kappa p}\leq 1 & \forall p\in V \label{eq:3}\\
    & \sum_{p=1}^{n} x_{1p}\leq 1 & \label{eq:4}\\
    & \sum_{p=1}^{n} x_{kp}\leq\sum_{p=1}^{n} x_{k-1,p} & \forall k\in\{2,\dots ,d+1\} \label{eq:5}\\
    & E_{1p}^\ell = x_{1p} & \forall \ell\in\{ s,t\},\ p\in V \label{eq:6}\\
    & \sum_{p=1}^{n} E_{kp}^\ell = 1 & \forall \ell\in\{ s,t\},\ k\in\{1,\dots ,d+1\} \label{eq:8}\\
    & E_{kp}^t + E_{kp}^s - x_{1p} \leq 1 & \forall k\in\{2,\dots ,d+1\},\ p\in V \label{eq:9}\\
    & E_{kp}^\ell \leq E_{k-1,p}^\ell + x_{kp} & \forall \ell\in\{ s,t\},\ k\in\{2,\dots ,d+1\},\ p\in V \label{eq:10}\\
    & x_{kp}\leq E_{kp}^t+E_{kp}^s & \forall k\in\{2,\dots ,d+1\},\ p\in V \label{eq:11}\\
    & E_{kp}^\ell\leq E_{k-1,p}^\ell+\sum_{k=1}^{n} a_{pk}\, E_{k-1,p}^\ell & \forall \ell\in\{ s,t\},\ k\in\{2,\dots ,d+1\},\ p\in V \label{eq:12}\\
    & x_{kp} \in\{0,1\} & \forall k\in\{1,\dots ,d+1\},\ p\in V \label{eq:13}\\
    & E_{kp}^\ell \in\{0,1\} & \forall \ell\in\{ s,t\},\ k\in\{1,\dots ,d+1\},\ p\in V \label{eq:14}
\end{align}\label{completeKnowledge}
\end{subequations}}

Constraints (\ref{eq:2}) ensure that only nodes that are available in iteration $k$ can be chosen. One node must not be chosen more than once (see (\ref{eq:3})) and there must not be more than one node chosen in iteration one (see (\ref{eq:4})). If there is one node that is chosen in iteration $k, k\geq 2$, then there must be a node that is chosen in the previous iteration $k-1$ (see (\ref{eq:5})). The combination of both, \eqref{eq:4} and \eqref{eq:5}, ensures that at most one point is chosen in each iteration. Both endpoints of the path in iteration $1$ are identical to the unique node that is chosen in iteration one (see (\ref{eq:6})). 
There must be exactly one start- and one endpoint of the path in every iteration (see (\ref{eq:8})). Both endpoints of a path must not be the same node in iteration $k\geq 2$ if one node is chosen in iteration $k$ (see (\ref{eq:9})). (Note that the constraint $E_{kp}^t + E_{kp}^s \leq 1$ is infeasible for an instance of the problem that has an optimal path of node length one.) 
A node can only become an endpoint in iteration $k$ if it was an endpoint in the previous iteration, or if it is chosen to extend the path in iteration $k$ (see (\ref{eq:10})). A node can only be chosen in iteration $k$, if it becomes an endpoint in iteration $k$  (see (\ref{eq:11})). A node may only become an endpoint in iteration $k$, if it was an endpoint in the previous iteration or if it is adjacent to an endpoint in the previous iteration (see (\ref{eq:12})), \(a_{pk}\) denotes the coefficient of the adjacency matrix. 


For a given instance of the problem, i.e., a given graph $G=(V,E)$ with $|V|=n$ and a given $(c,\Dupl)$-setting with associated availability matrix $D$, the above IP model  
has $O(d\cdot n)=O(\Dupl\cdot n^2)$ 
decision variables and 
$O(d\cdot n)=O(\Dupl\cdot n^2)$
constraints. Assuming that the number of duplicates $\Dupl$ is constant and does not grow with the problem size, both values are quadratic in the number of nodes in the graph.

The average computational time using CPLEX 12.10 to solve the IP is on average $0.6$~seconds when $n=22$, $\Dupl=3$ and $c=6$, which are the original parameters of the board game \TT. This computational time varies, if the parameters are different. For example, the IP solver needs on average $80$~seconds for one instance with $\Dupl=1$ and $c=6$ as parameter values. 
All numerical experiments were performed on a computer with an Intel Core i7-8700 CPU at 3.20\,GHz having 32\,GB RAM.


%% file: heuristics.tex
While the IP formulation given in Section~\ref{subsec:exsol} assumes complete knowledge of the order in which nodes become available during the game \TT, this is not the case in practice, i.e., when playing the game or when solving online path extension problems in general.  There is no complete knowledge in reality, 
and path extensions have to be chosen without the information on the nodes that become available in the next iteration(s). While being feasible in general, the solution of the IP formulation is an ideal solution that is generally hard to achieve in an online setting, but that certainly yields an upper bound on the best possible TT-path. In this section, we suggest path initialization heuristics (Iteration~1 of OPEP, selection of the first node) and associated path extension heuristics (Iterations~2 and following of OPEP, selection of consecutive nodes) to generate TT-paths in an online setting, i.e., when the initial tuple is known in iteration~1 and one further node becomes available in each consecutive iteration. Note that such path extension heuristics are online algorithms in the sense that they iteratively make decisions based on the partial knowledge available in a given iteration.

Throughout this section, we assume that a simple graph $G=(V,E)$ is given together with a feasible $(c,\Dupl)$-setting consisting of an initial set $\Vinit$ (or equivalently an initial tuple) and an extension tuple $(j_1,\dots ,j_d)$.

As mentioned above, we distinguish \emph{initialization heuristics} and \emph{extension heuristics}. An initialization heuristics initializes a TT-path by one initial node $\startp\in \Vinit$. The resulting initial path $P_0=(\startp)$ has node length one and is always feasible, i.e., it is always a TT-path. Both end nodes of this initial path are identical, i.e., $\startp=p_1=p_{\ell}$ in this particular case.

An extension heuristic assumes that, in iteration~$k$ of the OPEP with $k\geq 2$, a TT-path $P_{k-1}=(p_1,\dots,p_{k-1})$ has been computed during the preceding iterations, and that the set of available nodes in the current iteration is given by $\Vinit\cap\{j_1,\dots,j_{k} \}$. Then the extension heuristic suggests a strategy for the selection of the next node that is adjacent to $p_1$ or to $p_{k-1}$ and that is different from all nodes in $P_{k-1}$.

The following concepts will be useful for the formulation of initialization and extension heuristics.


\begin{definition}[Available tentacles]
    Let $G=(V,E)$ be a graph, let $P=(p_1,\dots,p_k)$ be a TT-path of node length $k\geq 1$, and let $\tentacleP$ be the set of tentacles of $P$ in $G$. Then 
    $\adVP\coloneqq\tentacleP\cap( \Vinit\cup\lbrace j_1,\dots ,j_{k}\})$ is the set of \emph{available tentacles} w.r.t.\ $P$, i.e., the set of nodes that are tentacles of $P$ and that are available in iteration $k$ of the OPEP.
\end{definition}
Note that the set of available tentacles comprises exactly those nodes that can be selected for a path extension in iteration $k$ of an extension heuristic. As a consequence, a further path extension is only possible if $\adVP\neq\emptyset$; otherwise the OPEP terminates.

\begin{definition}[Path extension]
    Let $G=(V,E)$ be a graph and let $P=(p_1,\dots,p_k)$ be a TT-path. A \textit{feasible extension} of $P$ is given by a tuple $(\extp,\pos)\in V\times\lbrace 0,1\rbrace$ such that  $\extp\in\adVP$.  Moreover, $\pos=1$ (indicating that $\extp$ is appended to $p_1$) requires that $(\extp,p_1)\in E$ and $\pos=0$ (indicating that $\extp$ is appended to $p_k$) requires that $(\extp,p_k)\in E$.
   The extended path is then given by $P\oplus\extp\coloneqq (p_1,\dots,p_{k},\extp)$ if $\pos=0$, and by $\extp\oplus P\coloneqq (\extp,p_1,\dots,p_{k})$ if $\pos=1$.  
\end{definition}



\subsubsection*{Random path initialization and random path extension}\label{subsec:random}

The simplest heuristic to initialize and to extend a (partial) TT-path is to always select a random node from the set of available nodes or the set of available tentacles, respectively.

\begin{algorithm}
	\caption{Start:Random}
	\label{heu: S1}
	\SetAlgoLined
\SetKwInOut{KwIn}{Input}
\SetKwInOut{KwOut}{Output}

\KwIn{Initial tuple $(i_1,\dots,i_c)$}
\KwOut{A node to initialize the TT-path: $P=(\startp)$}
 Compute a random integer $r$ between $1$ and $c$\;
 $\startp\coloneqq  i_r$.
\end{algorithm}


\begin{algorithm}
	\caption{Extension:Random}
	\label{heu: P1}
	\SetAlgoLined
\SetKwInOut{KwIn}{Input}
\SetKwInOut{KwOut}{Output}

\KwIn{Graph $G=(V,E)$, TT-path $P=(p_1,\dots,p_k)$, initial tuple $(i_1,\dots,i_c)$ and the first $k$ nodes of the extension tuple $(j_1,\dots ,j_{k})$}
\KwOut{A feasible extension $(\extp,\pos)$ or STOP}
    \uIf{$\adVP\neq\emptyset$}{
        Choose a random node $v_{ext}\in\adVP$\;
        \uIf{$(p_1,\extp)\in E \text{ and } (p_{k},\extp)\in E$}{
            Randomly choose $\pos\in\lbrace 0,1\rbrace$
        }
        \uElseIf{$(p_1,\extp)\in E$}{
            $\pos\coloneqq 1$
        }
        \Else{
            $\pos\coloneqq 0$
        }
    }
    \Else{
        STOP, there is no feasible extension.
    }
\end{algorithm}

\subsubsection*{Maximum degree heuristics}

When facing several choices in path initialization and path extension, we may prefer nodes with a large degree over nodes with a small degree. Indeed, the degree of a node may be interpreted as a (simple) indicator for the potential number of choices in consecutive iterations of OPEP. Moreover, this is a selection criterion that can be easily and efficiently evaluated.


\begin{algorithm}
	\caption{Start:Degree}
	\label{heu: S2}
	\SetAlgoLined
\SetKwInOut{KwIn}{Input}
\SetKwInOut{KwOut}{Output}

\KwIn{Graph $G=(V,E)$ and initial set $\Vinit$}
\KwOut{A node to initialize the TT-path: $P=(\startp)$}
 Compute $\delta (v)$ for every node $v\in\Vinit$\;
 $V_{\max}\coloneqq \argmax \lbrace \delta (v)\ :\ v\in \Vinit\rbrace$\;
 Choose a random node $\startp\in V_{\max}$.
\end{algorithm}



\begin{algorithm}
	\caption{Extension:Degree}
	\label{heu: P2}
	\SetAlgoLined
\SetKwInOut{KwIn}{Input}
\SetKwInOut{KwOut}{Output}

\KwIn{Graph $G=(V,E)$, TT-path $P=(p_1,\dots,p_k)$, initial tuple $(i_1,\dots,i_c)$ and the first $k$ nodes of the extension tuple $(j_1,\dots ,j_{k})$}
\KwOut{A feasible extension $(\extp,\pos)$ or STOP}
    \uIf{$\adVP\neq\emptyset$}{
        Compute $\delta (v)$ for every node $v\in\adVP$\;
        $V_{\max}\coloneqq \argmax\lbrace \delta (v)\ :\ v\in\adVP\rbrace$\;
        Choose a random node $\extp\in V_{\max}$\;
        \uIf{$(p_1,\extp)\in E \text{ and } (p_{k},\extp)\in E$}{
            Randomly choose $\pos\in\lbrace 0,1\rbrace$
        }
        \uElseIf{$(p_1,\extp)\in E$}{
            $\pos\coloneqq 1$
        }
        \Else{
            $\pos\coloneqq 0$
        }
    }
    \Else{
        STOP, there is no feasible extension.
    }
\end{algorithm}

\subsubsection*{Maximum tentacles heuristic}

A bit more subtle is the idea of using the number of tentacles of the extended path as a criterion for the next path extension.
Note that it does not make a difference whether the number of tentacles or the degree is used to select a starting node, see Lemma~\ref{lemma:bounds}. Thus, no new initialization heuristic is formulated for this case.

The number of tentacles of the extended path takes into account that not every node that is adjacent to a feasible extension of the current path is a feasible candidate for a consecutive extension (because it may, for example already be part of the path). Moreover, nodes are not counted twice if they are adjacent to both ends of the path.


Algorithm~\ref{heu: P3} summarizes the extension heuristic that uses tentacles as a selection criterion.

\begin{algorithm}
	\caption{Extension:Tentacles}
	\label{heu: P3}
	\SetAlgoLined
\SetKwInOut{KwIn}{Input}
\SetKwInOut{KwOut}{Output}

\KwIn{Graph $G=(V,E)$, TT-path $P=(p_1,\dots,p_k)$, initial tuple $(i_1,\dots,i_c)$ and the first $k$ node of the extension tuple $(j_1,\dots ,j_{k})$}
\KwOut{A feasible extension $(\extp,\pos)$ or STOP}
    \uIf{$\adVP\neq\emptyset$}{
        $\delta_{\max}\coloneqq  0$\;
        \For{$v\in\adVP$}{
            \If{$(p_1,v)\in E$}{
                \If{$\delta_{(v \oplus P)}>\delta_{\max}$}{
                    $\delta_{\max}\coloneqq \delta_{(v \oplus P)}$\;
                    $\extp\coloneqq  v$\;
                    $\pos\coloneqq  1$
                }
            }
            \If{$(p_{k},v)\in E$}{
                \If{$\delta_{(P \oplus v)}>\delta_{\max}$}{
                    $\delta_{\max}\coloneqq \delta_{(P \oplus v)}$\;
                    $\extp\coloneqq  v$\;
                    $\pos\coloneqq  0$
                }
            }
        }
    }
    \Else{
        STOP, there is no feasible extension.
    }
\end{algorithm}

\subsubsection{Connected components heuristics}

Whether there are feasible nodes to initialize or extend a path that are adjacent to further feasible extensions is an important information that can be used to potentially find better solutions, i.e., longer TT-paths. To formulate heuristics using this information, we need the following definition of a connected component.

\begin{definition}[Maximal connected component]\label{def: concom}
    Let $G=(V,E)$ be a graph and let $S\subseteq V$. 
    A set $Z\subseteq S$ is called a \emph{maximal connected component} of $S$  if $Z$ is connected in $S$ (i.e., every pair of nodes $a,b\in Z$ is connected by a path that lies completely in $S$) and if $Z$ is maximal with this property (i.e., there is no connected subset $Z'$ of $S$ with $Z\subsetneq Z'$).
    The set of all connected components of $S$ is denoted by $\mathcal{Z} (S)$.
\end{definition}


Note that a maximal connected component $Z$ of a subset $S\subseteq V$ may also be a singleton, i.e., it may consist of only one node.

Definition~\ref{def: concom} motivates the consideration of maximum connected components into which a current TT-path may be extended, i.e., maximal connected components with the largest possible number of nodes. If such a maximum connected component contains two or more nodes, we first select a node with minimal degree such that the nodes with larger degree are kept for later iterations. In the context of the initialization heuristic, we first consider $S\coloneqq\Vinit$.

\begin{algorithm}
	\caption{Start:Connected}
	\label{heu: S3}
	\SetAlgoLined
\SetKwInOut{KwIn}{Input}
\SetKwInOut{KwOut}{Output}

\KwIn{Graph $G=(V,E)$ and initial set $\Vinit$}
\KwOut{A node to initialize the TT-path: $P=(\startp)$}
 Compute $\mathcal{Z} (\Vinit)$ and set $z\coloneqq \max\lbrace |Z|\ :\ Z\in \mathcal{Z}(\Vinit) \rbrace$\;
    \uIf{$z = 1$}{
        Use Algorithm \ref{heu: S2}.
    }
    \Else{
        $Z_{\max}\coloneqq \bigcup \lbrace Z\ :\ |Z|=z \rbrace$\; 
        $V_{\min}\coloneqq \text{argmin}\lbrace \delta (v)\ :\ v\in Z_{\max}\rbrace$\;
        Choose a random node $\startp\in V_{\min}$.
    }
\end{algorithm}

When searching for a promising extension of a given TT-path $P=(p_1,\dots,p_k)$ in iteration $k\geq 2$, we look for maximum connected components that are in a sense adjacent to $P$. Towards this end, we consider maximum connected components solely in the set $V\setminus P\coloneqq V\setminus\{p_1,\dots,p_k\}$.


\begin{definition}[Maximal connected component excluding $P$]
    Let $G=(V,E)$ be a graph, let $P=(p_1,\dots,p_k)$ be a TT-path, and let $S\subseteq V$. A set $Z\subseteq (S\setminus P)$ is called a \textit{maximal connected component} of $S$ excluding $P$ if $Z$ is connected in $S\setminus P$ and if $Z$ is maximal with this property. The set of all maximal connected components of $S$ excluding $P$ is denoted by $\mathcal{Z}_P (S)$.
    Moreover, a set $Z\in\mathcal{Z}_P (S)$ is called \emph{adjacent} to $P$ if at least one element of $Z$ is adjacent to $p_1$ or to $p_{k}$. The set of all adjacent maximal connected components of $S$ excluding $P$ is denoted by $\adconP(S).$
\end{definition}

Note that an adjacent maximal connected component $Z$ of $S$ excluding $P$ consists only of nodes that are not part of the current TT-path $P$ and hence at least one node from $Z$ can be used to extend $P$. This motivates the formulation of an extension heuristic that aims to find extensions that are part of a largest possible adjacent connected component excluding $P$.

\begin{algorithm}
	\caption{Extension:Connected}
	\label{heu: P4}
	\SetAlgoLined
\SetKwInOut{KwIn}{Input}
\SetKwInOut{KwOut}{Output}

\KwIn{Graph $G=(V,E)$, TT-path $P=(p_1,\dots,p_k)$, initial tuple $(i_1,\dots,i_c)$ and the first $k$ nodes of the extension tuple $(j_1,\dots ,j_{k})$}
\KwOut{A feasible extension $(\extp,\pos)$ or STOP}
    \uIf{$\adVP\neq\emptyset$}{
        $S\coloneqq \lbrace i_1,\dots,i_c\rbrace\cup\lbrace j_1,\dots,j_{k}\rbrace$\;
        Compute $\adconP (S)$ and set $z\coloneqq \max\lbrace |Z|\ :\ Z\in \adconP (S) \rbrace$\;
        \uIf{$z =1$}{
            Use Algorithm \ref{heu: P3}
        }
        \Else{
             $Z_{\max}\coloneqq \bigcup \lbrace Z\ :\ |Z|=z \rbrace$\;
             $V_{\min}\coloneqq \text{argmin}\lbrace \delta (v)\ :\ v\in Z_{max}\cap \adVP\rbrace$\;
            Choose a random node $\extp\in V_{\min}$\;
            \uIf{$(p_1,\extp)\in E \text{ and } (p_{k},\extp)\in E$}{
                Randomly choose $\pos\in\lbrace 0,1\rbrace$
            }
            \uElseIf{$(p_1,\extp)\in E$}{
                Set $\pos\coloneqq 1$
            }
            \Else{
                Set $\pos\coloneqq 0$
            }
        }
    }
    \Else{
        STOP, there is no feasible extension.
    }
\end{algorithm}

\subsubsection*{Longest path heuristics}

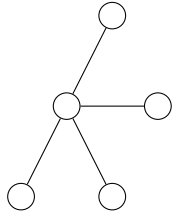
\begin{figure}
    \small
    \centering
    \begin{tikzpicture}[scale=0.6,every node/.style={draw=black,circle,font=\footnotesize,inner sep=1pt,minimum size=10pt}]
    \node (1) at (0,0) {};
    \node (2) at (-1,-2) {};
    \node (3) at (1,-2) {};
    \node (4) at (-2,-4) {};
    \node (5) at (0,-4) {};
    \draw[-] (1) to (2);
    \draw[-] (2) to (3);
    \draw[-] (2) to (4);
    \draw[-] (2) to (5);
    \end{tikzpicture}
    \caption{The node length of a path does not automatically grow with the cardinality of the connected component. The connected component shown here has cardinality $5$ but the path with maximum node length has only node length $3$.\label{Ex:CardinalityPath}}
\end{figure}

A further refinement of the connected components heuristics is obtained when, rather than the maximal connected components themselves, the node lengths of simple paths inside these components is considered. If there is a tie, the path with the larger number of tentacles should be used. Note that it is important to consider every  maximal (adjacent) connected component with cardinality greater or equal to $3$ (if there are any) to find a simple path with maximum node length in every situation. See Figure~\ref{Ex:CardinalityPath} for an illustration.

\begin{definition}[Node-longest paths extending $P$]
    Let $G=(V,E)$ be a graph, let $S\subseteq V$ and let $Z\in\mathcal{Z} (S)$. Then $\mathcal{P}_Z$ denotes the set of all \emph{node-longest paths} in $Z$. 
    Moreover, when  $P=(p_1,\dots,p_k)$ is a TT-path and $Z\in\adconP (S)$, then 
    $\longP$ denotes the set of all \emph{node-longest paths} in $P\cup Z$ that contain $P$ as a subpath. 
\end{definition}

Given a graph $G=(V,E)$, a subset $S\subseteq V$, a maximal connected component $Z\in\mathcal{Z}(S)$, and a path $P\in\longP$, it is easy to see that the following bounds are satisfied:
    \begin{enumerate}
        \item If $|Z|=k\leq 3$, then  $\ell (P)=k$.
        \item If $|Z|=k>3$, then $3\leq\ell (P)\leq k$.
    \end{enumerate}
See again Figure~\ref{Ex:CardinalityPath} for an illustration.


\begin{algorithm}
	\caption{Start:LongestPath}
	\label{heu: S4}
	\SetAlgoLined
\SetKwInOut{KwIn}{Input}
\SetKwInOut{KwOut}{Output}

\KwIn{Graph $G=(V,E)$ and initial set $\Vinit$}
\KwOut{A node to initialize the TT-path: $P=(\startp)$}
 Compute $\mathcal{Z} (\Vinit)$ and set $z\coloneqq \max\lbrace |Z|\ :\ Z\in \mathcal{Z}(\Vinit) \rbrace$\;
    \uIf{$z = 1$}{
        Use Algorithm \ref{heu: S2}.
        }
    \Else{
    $\kappa\coloneqq\min\{z,3\}$
    }
     $\ell_{\max}\coloneqq 0,\ \delta_{\max}\coloneqq 0,\ Z^*\coloneqq \emptyset$\;
    \For{$Z\in\mathcal{Z} (\Vinit)$ with $|Z|\geq \kappa$}{
        $P_Z^\delta\coloneqq \lbrace \bar{P}\in\mathcal{P}_Z\ :\ \delta_{\bar{P}}\geq\delta_{P'}\ \forall P'\in\mathcal{P}_Z\rbrace$\;
        Choose a random path $P\in P_Z^\delta$\;
        \If{$(\ \ell(P)>\ell_{\max}\ )\ \text{or} (\ \ell(P)=\ell_{\max}\ \text{and}\ \delta_P>\delta_{\max}\ )$}{
            $\ell_{\max}\coloneqq \ell(P),\ \delta_{\max}\coloneqq \delta_P,\ Z^*\coloneqq Z$\;
        }
    }
    $V_{\min}\coloneqq \text{argmin}\lbrace \delta (v)\ :\ v\in P_{Z^*}^\delta\rbrace$\;
    Choose a random node $\startp\in V_{\min}$.
\end{algorithm}

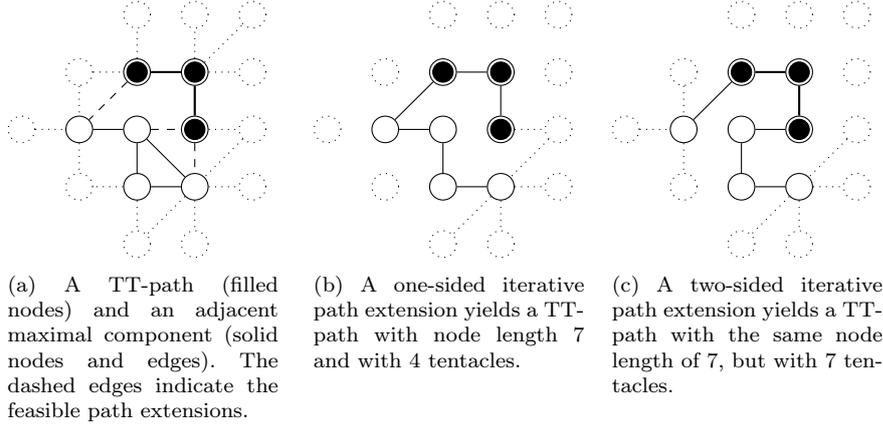
\begin{figure}
    \small\centering
    \hspace{\fill}
    \begin{subfigure}[t]{0.3\textwidth}
    \begin{tikzpicture}[scale=0.38,every node/.style={draw=black,circle,font=\footnotesize,inner sep=1pt,minimum size=10pt}]
    \node[dotted] (1) at (2,0) {};
    \node (2) at (4,0) {};
    \node (3) at (6,0) {};
    \node[dotted] (4) at (8,0) {};
    \node[dotted] (5) at (0,-2) {};
    \node (6) at (2,-2) {};
    \node (7) at (4,-2) {};
    \node (8) at (6,-2) {};
    \node[dotted] (9) at (8,-2) {};
    \node (10) at (4,-4) {};
    \node (11) at (6,-4) {};
    \node[dotted] (12) at (8,-4) {};
    \node[dotted] (13) at (2,-4) {};
    \node[dotted] (14) at (4,-6) {};
    \node[dotted] (15) at (6,-6) {};
    \node[dotted] (16) at (4,2) {};
    \node[dotted] (17) at (6,2) {};
    \node[dotted] (18) at (8,2) {};
    \filldraw[black] (2) circle (10pt);
    \filldraw[black] (3) circle (10pt);
    \filldraw[black] (8) circle (10pt);
    \draw[-,thick] (2) to (3);
    \draw[-,thick] (3) to (8);
    \draw[-] (6) to (7);
    \draw[-] (7) to (10);
    \draw[-] (7) to (11);
    \draw[-] (10) to (11);
    \draw[dashed] (2) to (6);
    \draw[dashed] (8) to (11);
    \draw[dashed] (7) to (8);
    \draw[dotted] (1) to (2);
    \draw[dotted] (1) to (6);
    \draw[dotted] (2) to (16);
    \draw[dotted] (3) to (4);
    \draw[dotted] (3) to (17);
    \draw[dotted] (3) to (18);
    \draw[dotted] (5) to (6);
    \draw[dotted] (6) to (13);
    \draw[dotted] (8) to (9);
    \draw[dotted] (9) to (11);
    \draw[dotted] (10) to (13);
    \draw[dotted] (10) to (14);
    \draw[dotted] (11) to (12);
    \draw[dotted] (11) to (14);
    \draw[dotted] (11) to (15);
    \end{tikzpicture}
    \label{Ex:LongstSubpathA}
    \caption{A TT-path (filled nodes) and an adjacent maximal component (solid nodes and edges). The dashed edges indicate the feasible path extensions.}
    \end{subfigure}
    \hspace{\fill}
    \begin{subfigure}[t]{0.3\textwidth}
    \begin{tikzpicture}[scale=0.38,every node/.style={draw=black,circle,font=\footnotesize,inner sep=1pt,minimum size=10pt}]
    \node[dotted] (1) at (2,0) {};
    \node (2) at (4,0) {};
    \node (3) at (6,0) {};
    \node[dotted] (4) at (8,0) {};
    \node[dotted] (5) at (0,-2) {};
    \node (6) at (2,-2) {};
    \node (7) at (4,-2) {};
    \node (8) at (6,-2) {};
    \node[dotted] (9) at (8,-2) {};
    \node (10) at (4,-4) {};
    \node (11) at (6,-4) {};
    \node[dotted] (12) at (8,-4) {};
    \node[dotted] (13) at (2,-4) {};
    \node[dotted] (14) at (4,-6) {};
    \node[dotted] (15) at (6,-6) {};
    \node[dotted] (16) at (4,2) {};
    \node[dotted] (17) at (6,2) {};
    \node[dotted] (18) at (8,2) {};
    \filldraw[black] (2) circle (10pt);
    \filldraw[black] (3) circle (10pt);
    \filldraw[black] (8) circle (10pt);
    \draw[-] (2) to (3);
    \draw[-] (3) to (8);
    \draw[-] (6) to (7);
    \draw[-] (7) to (10);
    \draw[-] (10) to (11);
    \draw[-] (2) to (6);
    \draw[dotted] (8) to (9);
    \draw[dotted] (9) to (11);
    \draw[dotted] (11) to (12);
    \draw[dotted] (11) to (14);
    \draw[dotted] (11) to (15);
    \end{tikzpicture}
    \label{Ex:LongstSubpathB}
    \caption{A one-sided iterative  path extension yields a TT-path with node length $7$ and with $4$ tentacles.}
    \end{subfigure}
    \hfill
    \begin{subfigure}[t]{0.3\textwidth}
    \begin{tikzpicture}[scale=0.38,every node/.style={draw=black,circle,font=\footnotesize,inner sep=1pt,minimum size=10pt}]
    \node[dotted] (1) at (2,0) {};
    \node (2) at (4,0) {};
    \node (3) at (6,0) {};
    \node[dotted] (4) at (8,0) {};
    \node[dotted] (5) at (0,-2) {};
    \node (6) at (2,-2) {};
    \node (7) at (4,-2) {};
    \node (8) at (6,-2) {};
    \node[dotted] (9) at (8,-2) {};
    \node (10) at (4,-4) {};
    \node (11) at (6,-4) {};
    \node[dotted] (12) at (8,-4) {};
    \node[dotted] (13) at (2,-4) {};
    \node[dotted] (14) at (4,-6) {};
    \node[dotted] (15) at (6,-6) {};
    \node[dotted] (16) at (4,2) {};
    \node[dotted] (17) at (6,2) {};
    \node[dotted] (18) at (8,2) {};
    \filldraw[black] (2) circle (10pt);
    \filldraw[black] (3) circle (10pt);
    \filldraw[black] (8) circle (10pt);
    \draw[-,thick] (2) to (3);
    \draw[-,thick] (3) to (8);
    \draw[-] (7) to (10);
    \draw[-] (10) to (11);
    \draw[-] (2) to (6);
    \draw[-] (7) to (8);
    \draw[dotted] (1) to (6);
    \draw[dotted] (5) to (6);
    \draw[dotted] (6) to (13);
    \draw[dotted] (9) to (11);
    \draw[dotted] (11) to (12);
    \draw[dotted] (11) to (14);
    \draw[dotted] (11) to (15);
    \end{tikzpicture}
    \label{Ex:LongstSubpathC}
    \caption{A two-sided iterative path extension yields a TT-path with the same node length of $7$, but with $7$ tentacles.}
    \end{subfigure}
    \caption{Illustration of the importance of allowing a given TT-path $P$ to be extended on both ends while computing a node-longest possible path inside a maximal connected component. Filled nodes represent a current TT-path $P$. Dotted nodes are not part of the connected component nor the given path. \label{Ex:LongestSubpath}}
\end{figure}

The same idea can be used to formulate an extension heuristic, see Algorithm~\ref{heu: P5}. Note that it is important to allow the extension of the path on both ends to find the best solution inside one connected component. Of course this is only important if one evaluated connected component is adjacent to both endpoints of the given path. See Figure~\ref{Ex:LongestSubpath} for an illustration.

\begin{algorithm}
	\caption{Extension:LongestPath}
	\label{heu: P5}
	\SetAlgoLined
\SetKwInOut{KwIn}{Input}
\SetKwInOut{KwOut}{Output}

\KwIn{Graph $G=(V,E)$, TT-path $P=(p_1,\dots,p_k)$, initial tuple $(i_1,\dots,i_c)$ and first $k$ nodes of the extension tuple $(j_1,\dots ,j_{k})$}
\KwOut{A feasible extension $(\extp,\pos)$ or STOP}
    \uIf{$\adVP\neq\emptyset$}{
        $S\coloneqq \lbrace i_1,\dots,i_c\rbrace\cup\lbrace j_1,\dots ,j_{k}\rbrace$\;
        Compute $\adconP (S)$ and set $z\coloneqq \max\lbrace |Z|\ :\ Z\in \adconP (S) \rbrace$\;
        \uIf{$z =1$}{
            Use Algorithm \ref{heu: P3}
        }
        \Else{
        $\kappa\coloneqq \min\{z,3\}$
        }
        $\ell_{\max}=0,\ \delta_{\max},\ P^*=\emptyset$\;
        \For{$Z\in\adconP (S)$ with $|Z|\geq \kappa$}{
            $P_Z^\delta\coloneqq \lbrace P'\in\longP\ :\ \delta_{P'}\geq\delta_{\bar{P}}\ \forall \bar{P}\in\longP\rbrace$\;
            Choose a random $P'\in P_Z^\delta$\;
            \uIf{$\ell (P')>\ell_{\max}$}{
                $P^*\coloneqq P'$, $\ell_{\max}\coloneqq \ell (P')$, $\delta_{\max}\coloneqq \delta_{P'}$
            }
            \ElseIf{$\ell (P')=\ell_{\max}\text{ and }\delta_{P'}>\delta_{\max}$}{
                    $P^*\coloneqq P'$, $\ell_{\max}\coloneqq \ell (P')$, $\delta_{\max}\coloneqq \delta_{P'}$
            }
            Choose $\extp\in P^*$ such that it is adjacent to $p_1$ ($\pos =1$) or $p_k$ ($\pos =0$) in the path $P^*$ and the selection maximizes the number of tentacles of the resulting path\;
        }
    }
    \Else{
        STOP, there is no feasible extension.
    }
\end{algorithm}

\subsubsection*{Illustration and comparison}

\begin{figure}
    \small
    \centering
    \hspace{\fill}
    \begin{subfigure}[t]{0.45\textwidth}
    \begin{tikzpicture}[scale=0.6,every node/.style={draw=black,circle,font=\footnotesize,inner sep=1pt,minimum size=10pt}]
    \node[dotted] (-2) at (-2,4) {};
    \node[dotted] (-1) at (-2,2) {};
    \node[dotted] (0) at (-2,0) {};
    \node[dotted] (1) at (0,0) {};
    \node[dotted] (2) at (2,0) {};
    \node (3) at (4,0) {};
    \node (4) at (6,0) {};
    \node[dotted] (5) at (0,2) {};
    \node[dashed] (6) at (2,2) {};
    \node (7) at (4,2) {};
    \node[dotted] (8) at (6,2) {};
    \node (9) at (0,4) {};
    \node (10) at (2,4) {};
    \node (11) at (4,4) {};
    \node[dotted] (12) at (6,4) {};
    \filldraw[black] (2) circle (10pt);
    \filldraw[black] (5) circle (10pt);
    \filldraw[black] (10) circle (10pt);
    \draw[dotted] (-2) to (9);
    \draw[dotted] (-1) to (-2);
    \draw[dotted] (-1) to (5);
    \draw[dotted] (0) to (-1);
    \draw[dotted] (0) to (1);
    \draw[dotted] (0) to (5);
    \draw[dotted] (2) to (3);
    \draw[-] (2) to (5);
    \draw[dotted] (2) to (6);
    \draw[dotted] (3) to (7);
    \draw[dotted] (4) to (7);
    \draw[dotted] (5) to (6);
    \draw[dotted] (5) to (9);
    \draw[-] (5) to (10);
    \draw[dotted] (6) to (7);
    \draw[dashed] (6) to (10);
    \draw[dotted] (7) to (8);
    \draw[dotted] (7) to (11);
    \draw[dotted] (8) to (12);
    \draw[dotted] (9) to (10);
    \draw[dotted] (10) to (11);
    \draw[dotted] (11) to (12);
    \end{tikzpicture}
    \caption{Extension:Degree selects an available node with maximum degree}
    \end{subfigure}
    \hspace{\fill}
    \begin{subfigure}[t]{0.45\textwidth}\centering
    \begin{tikzpicture}[scale=0.6,every node/.style={draw=black,circle,font=\footnotesize,inner sep=1pt,minimum size=10pt}]
    \node[dotted] (-2) at (-2,4) {};
    \node[dotted] (-1) at (-2,2) {};
    \node[dotted] (0) at (-2,0) {};
    \node[dotted] (1) at (0,0) {};
    \node[dotted] (2) at (2,0) {};
    \node[dashed] (3) at (4,0) {};
    \node (4) at (6,0) {};
    \node[dotted] (5) at (0,2) {};
    \node (6) at (2,2) {};
    \node (7) at (4,2) {};
    \node[dotted] (8) at (6,2) {};
    \node (9) at (0,4) {};
    \node (10) at (2,4) {};
    \node (11) at (4,4) {};
    \node[dotted] (12) at (6,4) {};
    \filldraw[black] (2) circle (10pt);
    \filldraw[black] (5) circle (10pt);
    \filldraw[black] (10) circle (10pt);
    \draw[dotted] (-2) to (9);
    \draw[dotted] (-1) to (-2);
    \draw[dotted] (-1) to (5);
    \draw[dotted] (0) to (-1);
    \draw[dotted] (0) to (1);
    \draw[dotted] (0) to (5);
    \draw[dashed] (2) to (3);
    \draw[-] (2) to (5);
    \draw[dotted] (2) to (6);
    \draw[dotted] (3) to (7);
    \draw[dotted] (4) to (7);
    \draw[dotted] (5) to (6);
    \draw[dotted] (5) to (9);
    \draw[-] (5) to (10);
    \draw[dotted] (6) to (7);
    \draw[dotted] (6) to (10);
    \draw[dotted] (7) to (8);
    \draw[dotted] (7) to (11);
    \draw[dotted] (8) to (12);
    \draw[dotted] (9) to (10);
    \draw[dotted] (10) to (11);
    \draw[dotted] (11) to (12);
    \end{tikzpicture}
    \caption{Extension:Tentacle selects an available node that maximizes the number of tentacles of the extended path}
    \end{subfigure}
    \hspace{\fill}
    \\[4ex]
    \hspace{\fill}
    \begin{subfigure}[t]{0.45\textwidth}
    \begin{tikzpicture}[scale=0.6,every node/.style={draw=black,circle,font=\footnotesize,inner sep=1pt,minimum size=10pt}]
    \node[dotted] (-2) at (-2,4) {};
    \node[dotted] (-1) at (-2,2) {};
    \node[dotted] (0) at (-2,0) {};
    \node[dotted] (1) at (0,0) {};
    \node[dotted] (2) at (2,0) {};
    \node[dashed] (3) at (4,0) {};
    \node (4) at (6,0) {};
    \node[dotted] (5) at (0,2) {};
    \node (6) at (2,2) {};
    \node (7) at (4,2) {};
    \node[dotted] (8) at (6,2) {};
    \node (9) at (0,4) {};
    \node (10) at (2,4) {};
    \node (11) at (4,4) {};
    \node[dotted] (12) at (6,4) {};
    \filldraw[black] (2) circle (10pt);
    \filldraw[black] (5) circle (10pt);
    \filldraw[black] (10) circle (10pt);
    \draw[dotted] (-2) to (9);
    \draw[dotted] (-1) to (-2);
    \draw[dotted] (-1) to (5);
    \draw[dotted] (0) to (-1);
    \draw[dotted] (0) to (1);
    \draw[dotted] (0) to (5);
    \draw[dashed] (2) to (3);
    \draw[-] (2) to (5);
    \draw[dotted] (2) to (6);
    \draw[-] (3) to (7);
    \draw[-] (4) to (7);
    \draw[dotted] (5) to (6);
    \draw[dotted] (5) to (9);
    \draw[-] (5) to (10);
    \draw[-] (6) to (7);
    \draw[dotted] (6) to (10);
    \draw[dotted] (7) to (8);
    \draw[-] (7) to (11);
    \draw[dotted] (8) to (12);
    \draw[dotted] (9) to (10);
    \draw[dotted] (10) to (11);
    \draw[dotted] (11) to (12);
    \end{tikzpicture}
    \caption{Extension:Connected selects an adjacent node from a maximum adjacent component}
    \end{subfigure}
    \hspace{\fill}
    \begin{subfigure}[t]{0.45\textwidth}\centering
    \begin{tikzpicture}[scale=0.6,every node/.style={draw=black,circle,font=\footnotesize,inner sep=1pt,minimum size=10pt}]
    \node[dotted] (-2) at (-2,4) {};
    \node[dotted] (-1) at (-2,2) {};
    \node[dotted] (0) at (-2,0) {};
    \node[dotted] (1) at (0,0) {};
    \node[dotted] (2) at (2,0) {};
    \node[dashed] (3) at (4,0) {};
    \node (4) at (6,0) {};
    \node[dotted] (5) at (0,2) {};
    \node (6) at (2,2) {};
    \node (7) at (4,2) {};
    \node[dotted] (8) at (6,2) {};
    \node (9) at (0,4) {};
    \node (10) at (2,4) {};
    \node (11) at (4,4) {};
    \node[dotted] (12) at (6,4) {};
    \filldraw[black] (2) circle (10pt);
    \filldraw[black] (5) circle (10pt);
    \filldraw[black] (10) circle (10pt);
    \draw[dotted] (-2) to (9);
    \draw[dotted] (-1) to (-2);
    \draw[dotted] (-1) to (5);
    \draw[dotted] (0) to (-1);
    \draw[dotted] (0) to (1);
    \draw[dotted] (0) to (5);
    \draw[dashed] (2) to (3);
    \draw[-] (2) to (5);
    \draw[dotted] (2) to (6);
    \draw[dashed] (3) to (7);
    \draw[-] (4) to (7);
    \draw[dotted] (5) to (6);
    \draw[dotted] (5) to (9);
    \draw[-] (5) to (10);
    \draw[-] (6) to (7);
    \draw[dashed] (6) to (10);
    \draw[dotted] (7) to (8);
    \draw[dashed] (7) to (11);
    \draw[dotted] (8) to (12);
    \draw[dotted] (9) to (10);
    \draw[dotted] (10) to (11);
    \draw[dotted] (11) to (12);
    \end{tikzpicture}
    \caption{Extension:LongestPath computes a maximum possible (two-sided) path extension as the basis for the selection}
    \end{subfigure}
    \hspace{\fill}
    \caption{Illustration of different extension heuristics at an example with a current TT-path of node length $3$ (filled nodes) and with $6$ available nodes (solid and dashed nodes). Dotted nodes are not available. The selected extension is dashed and the planned path is highlighted by dashed edges.\label{Ex:ExtHeu}}
\end{figure}
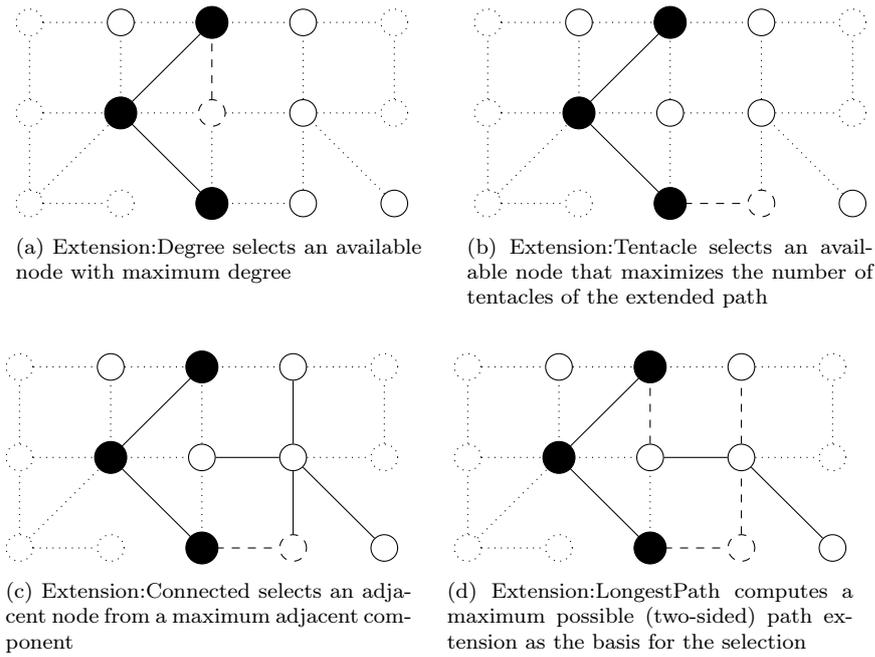

Figures~\ref{Ex:StartHeu} and \ref{Ex:ExtHeu} illustrate the selections made by the different initialization and extension heuristics introduced in the previous sections. It can be expected that the more information is used for the selection of the next node in the OPEP, the higher is the probability of generating longer TT-paths. While this is confirmed in general by the numerical studies presented in Section~\ref{sec:res}, it should be noted that for each of the heuristics it is possible to construct instances where the respective performance is very bad compared to the ideal solution. Moreover, a more involved heuristic naturally needs higher computational times, which may become a critical issue when larger instances are considered in an online setting.

\begin{figure}
    \small
    \centering
    \hspace{\fill}
    \begin{subfigure}{0.4\textwidth}
    \begin{tikzpicture}[scale=0.6,every node/.style={draw=black,circle,font=\footnotesize,inner sep=1pt,minimum size=10pt}]
    \node[dotted] (-2) at (-2,4) {};
    \node[dotted] (-1) at (-2,2) {};
    \node[dotted] (0) at (-2,0) {};
    \node[dotted] (1) at (0,0) {};
    \node[dotted] (2) at (2,0) {};
    \node (3) at (4,0) {};
    \node (4) at (6,0) {};
    \node[dashed] (5) at (0,2) {};
    \node[dotted] (6) at (2,2) {};
    \node (7) at (4,2) {};
    \node (8) at (6,2) {};
    \node[dotted] (9) at (0,4) {};
    \node[dotted] (10) at (2,4) {};
    \node[dotted] (11) at (4,4) {};
    \node (12) at (6,4) {};
    \draw[dotted] (-2) to (9);
    \draw[dotted] (-1) to (-2);
    \draw[dotted] (-1) to (5);
    \draw[dotted] (0) to (-1);
    \draw[dotted] (0) to (1);
    \draw[dotted] (0) to (5);
    \draw[dotted] (2) to (3);
    \draw[dotted] (2) to (5);
    \draw[dotted] (2) to (6);
    \draw[dotted] (3) to (7);
    \draw[dotted] (4) to (7);
    \draw[dotted] (5) to (6);
    \draw[dotted] (5) to (9);
    \draw[dotted] (5) to (10);
    \draw[dotted] (6) to (7);
    \draw[dotted] (6) to (10);
    \draw[dotted] (7) to (8);
    \draw[dotted] (7) to (11);
    \draw[dotted] (8) to (12);
    \draw[dotted] (9) to (10);
    \draw[dotted] (10) to (11);
    \draw[dotted] (11) to (12);
    \end{tikzpicture}
    \caption{Start:Degree selects the node with the largest degree}
    \end{subfigure}
    \hspace{\fill}
    \begin{subfigure}{0.4\textwidth}
    \begin{tikzpicture}[scale=0.6,every node/.style={draw=black,circle,font=\footnotesize,inner sep=1pt,minimum size=10pt}]
    \node[dotted] (-2) at (-2,4) {};
    \node[dotted] (-1) at (-2,2) {};
    \node[dotted] (0) at (-2,0) {};
    \node[dotted] (1) at (0,0) {};
    \node[dotted] (2) at (2,0) {};
    \node (3) at (4,0) {};
    \node[dashed] (4) at (6,0) {};
    \node (5) at (0,2) {};
    \node[dotted] (6) at (2,2) {};
    \node (7) at (4,2) {};
    \node (8) at (6,2) {};
    \node[dotted] (9) at (0,4) {};
    \node[dotted] (10) at (2,4) {};
    \node[dotted] (11) at (4,4) {};
    \node (12) at (6,4) {};
    \draw[dotted] (-2) to (9);
    \draw[dotted] (-1) to (-2);
    \draw[dotted] (-1) to (5);
    \draw[dotted] (0) to (-1);
    \draw[dotted] (0) to (1);
    \draw[dotted] (0) to (5);
    \draw[dotted] (2) to (3);
    \draw[dotted] (2) to (5);
    \draw[dotted] (2) to (6);
    \draw[-] (3) to (7);
    \draw[-] (4) to (7);
    \draw[dotted] (5) to (6);
    \draw[dotted] (5) to (9);
    \draw[dotted] (5) to (10);
    \draw[dotted] (6) to (7);
    \draw[dotted] (6) to (10);
    \draw[-] (7) to (8);
    \draw[dotted] (7) to (11);
    \draw[-] (8) to (12);
    \draw[dotted] (9) to (10);
    \draw[dotted] (10) to (11);
    \draw[dotted] (11) to (12);
    \end{tikzpicture}
    \caption{Start:Connected selects a node with minimum degree in a maximum connected component in $\Vinit$}
    \end{subfigure}
    \hspace{\fill}
    \\[3ex]
    \begin{subfigure}{0.4\textwidth}
    \begin{tikzpicture}[scale=0.6,every node/.style={draw=black,circle,font=\footnotesize,inner sep=1pt,minimum size=10pt}]
    \node[dotted] (-2) at (-2,4) {};
    \node[dotted] (-1) at (-2,2) {};
    \node[dotted] (0) at (-2,0) {};
    \node[dotted] (1) at (0,0) {};
    \node[dotted] (2) at (2,0) {};
    \node (3) at (4,0) {};
    \node (4) at (6,0) {};
    \node (5) at (0,2) {};
    \node[dotted] (6) at (2,2) {};
    \node (7) at (4,2) {};
    \node[dashed] (8) at (6,2) {};
    \node[dotted] (9) at (0,4) {};
    \node[dotted] (10) at (2,4) {};
    \node[dotted] (11) at (4,4) {};
    \node (12) at (6,4) {};
    \draw[dotted] (-2) to (9);
    \draw[dotted] (-1) to (-2);
    \draw[dotted] (-1) to (5);
    \draw[dotted] (0) to (-1);
    \draw[dotted] (0) to (1);
    \draw[dotted] (0) to (5);
    \draw[dotted] (2) to (3);
    \draw[dotted] (2) to (5);
    \draw[dotted] (2) to (6);
    \draw[dashed] (3) to (7);
    \draw[-] (4) to (7);
    \draw[dotted] (5) to (6);
    \draw[dotted] (5) to (9);
    \draw[dotted] (5) to (10);
    \draw[dotted] (6) to (7);
    \draw[dotted] (6) to (10);
    \draw[dashed] (7) to (8);
    \draw[dotted] (7) to (11);
    \draw[dashed] (8) to (12);
    \draw[dotted] (9) to (10);
    \draw[dotted] (10) to (11);
    \draw[dotted] (11) to (12);
    \end{tikzpicture}
    \caption{Start:LongestPath identifies the longest path in a maximum connected component in $\Vinit$}
    \end{subfigure}
    \caption{Illustration of different starting heuristics at an example with $|\Vinit|=6$. Dotted nodes are in $V\setminus\Vinit$ and  can thus not be chosen. The dashed node is selected by the respective starting heuristic.\label{Ex:StartHeu}}
\end{figure}
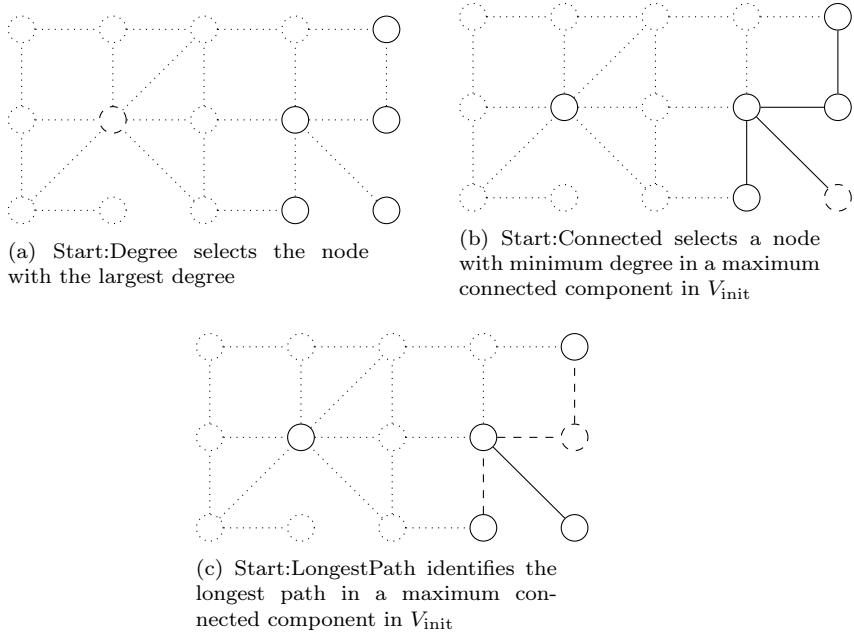

\begin{figure}
    \small
    \centering
    \begin{tikzpicture}[scale=0.6,every node/.style={draw=black,circle,font=\footnotesize,inner sep=1pt,minimum size=15pt}]
    \node (1) at (0,0) {9};
    \node (2) at (2,0) {1};
    \node (3) at (4,0) {10};
    \node (4) at (6,0) {2};
    \node (5) at (8,0) {1};
    \node (6) at (0,-2) {1};
    \node (7) at (2,-2) {16};
    \node (8) at (4,-2) {17};
    \node (9) at (6,-2) {11};
    \node (10) at (8,-2) {3};
    \node (11) at (-2,-4) {8};
    \node (12) at (0,-4) {15};
    \node (13) at (2,-4) {1};
    \node (14) at (4,-4) {14};
    \node (15) at (6,-4) {13};
    \node (16) at (8,-4) {12};
    \node (17) at (10,-4) {4};
    \node (18) at (-2,-6) {1};
    \node (19) at (0,-6) {7};
    \node (20) at (2,-6) {6};
    \node (21) at (4,-6) {5};
    \node (22) at (6,-6) {1};
    \draw[-] (1) to (2);
    \draw[-] (1) to (6);
    \draw[-] (1) to (7);
    \draw[-] (2) to (3);
    \draw[-] (2) to (7);
    \draw[-] (3) to (4);
    \draw[-] (3) to (7);
    \draw[-] (3) to (8);
    \draw[-] (3) to (9);
    \draw[-] (4) to (5);
    \draw[-] (4) to (9);
    \draw[-] (4) to (10);
    \draw[-] (6) to (7);
    \draw[-] (6) to (11);
    \draw[-] (7) to (8);
    \draw[-] (7) to (12);
    \draw[-] (7) to (13);
    \draw[-] (8) to (9);
    \draw[-] (8) to (13);
    \draw[-] (8) to (14);
    \draw[-] (8) to (15);
    \draw[-] (9) to (15);
    \draw[-] (9) to (16);
    \draw[-] (10) to (17);
    \draw[-] (11) to (12);
    \draw[-] (11) to (18);
    \draw[-] (11) to (19);
    \draw[-] (12) to (13);
    \draw[-] (12) to (19);
    \draw[-] (12) to (20);
    \draw[-] (13) to (14);
    \draw[-] (13) to (20);
    \draw[-] (14) to (15);
    \draw[-] (14) to (20);
    \draw[-] (14) to (21);
    \draw[-] (15) to (16);
    \draw[-] (15) to (21);
    \draw[-] (15) to (22);
    \draw[-] (16) to (17);
    \draw[-] (16) to (22);
    \draw[-] (17) to (22);
    \draw[-] (18) to (19);
    \draw[-] (19) to (20);
    \draw[-] (20) to (21);
    \draw[-] (21) to (22);
    \end{tikzpicture}
    \caption{Illustration of a worst-case scenario. The numbers in the nodes indicate the iteration in which the respective node becomes available. The initial tuple consists of six nodes (nodes with value $1$) which define an independent set in $G$. Nevertheless, the ideal solution (which can be, for example, computed by solving the IP formulation from Section~\ref{subsec:exsol}) has node length $22$. Every starting heuristic (neglecting Start:Random) will lead to a path of node length $1$. Start:Random has a $\frac{1}{6}$ chance to select the ``correct'' starting node.\label{Ex:WorstCase}}
\end{figure}
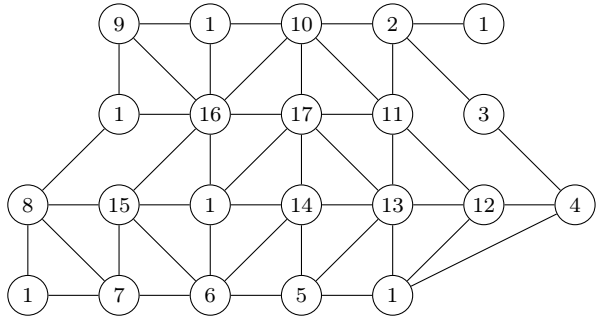

We provide an example of a scenario in which most initialization heuristics do not find a reasonable solution in Figure~\ref{Ex:WorstCase}. Indeed, whenever the initial set is an independent set and if, in addition, the first vertex of the extension tuple is adjacent to one of the nodes of the initial set, an initialization heuristic will only lead to an extendable path if the \emph{one} ``correct'' node is selected. This, however, can not be expected in general.

\subsection{Combinations of starting heuristics and extension heuristics}\label{subsubsec:combheu}

While the initialization heuristics and the extension heuristics come in pairs of methods based on similar ideas, other combinations of initialization and extension heuristics are of course possible. Moreover, every extension could be realized by using an iteration of a different extension heuristic. The general concept is summarized in Algorithm~\ref{heu:combined}.

\begin{algorithm}
	\caption{Combined Heuristic}
	\label{heu:combined}
	\SetAlgoLined
\SetKwInOut{KwIn}{Input}
\SetKwInOut{KwOut}{Output}

\KwIn{Graph $G=(V,E)$, initial tuple $(i_1,\dots,i_c)$ and extension tuple $(j_1,\dots ,j_d)$}
\KwOut{A TT-path $P$}
    Use an initialization heuristic to compute $\startp$\;
    $P=(\startp)$\;
    \Repeat{STOP}{
        Use an extension heuristic to compute $(\extp,\pos)$\;
        \If{STOP}{
            Return $P$.
        }
        \uIf{$\pos=1$}{
            $P=(\extp\oplus P)$.
        }\Else{
            $P=(P\oplus\extp)$.
        }
    }
    
\end{algorithm}

%% file: Results.tex
In this section we test, compare and discuss the different solution heuristics. In this context the solutions IP model are
used
as reference solutions to evaluate the quality of the solutions of the different heuristics. The comparison
between the optimal solutions of the IP model and the solutions of the heuristic approaches is not fair in two respects. First the IP model is based on complete knowledge of the card deck, while the heuristics have partial knowledge namely the currently available extensions in each iteration. Second, obviously construction heuristics are by orders of magnitude faster than IP solvers and can thus not be expected to achieve the solution quality.

\subsection{Variation of problem parameters}\label{subsec:normal_diffpar}
In our numerical test we evaluate the proposed heuristics on different types of instances. 
Problem parameters like the size of the considered graph, the size of the initial tuple or the number of duplicates effect both solution quality and running time of the heuristics.

\subsubsection*{Number of nodes and edges in $G$}\label{subsubsec:G}
The underlying graph has a significant impact on the run of the game and the performance of the heuristics. The first graph we use for our numerical tests is the original graph \G\ of the board game \TT\ which consists of $22$~nodes and $45$~edges, see Figure~\ref{Fig:Boardgame as graph}. Furthermore, we apply the proposed heuristics on an enlarged version of this graph \Gext\ (see Figure~\ref{fig:largergraph}), which consists of $32$ nodes and $66$ edges and has a similar connectivity structure as the original one. The enlarged graph contains roughly $50\%$ more nodes and edges than the original graph of the board game. 
The longest simple paths contain in both cases every node of the graph, which is an upper bound with respect to all feasible settings.

\subsubsection*{Different values of $c$}\label{subsubsec:c}
In the setup of the board game \TT\ the initial tuple has size $c=6$. We test our heuristics on the original board game graph \G\ with the values \(c\in\{4,\ldots,8\}\). On the extended graph \Gext\ with $32$ nodes we chose \(c\in\{7,\ldots,11\}\).
Note that every solution found for a fixed value of $c$ is also a feasible solution for $c=c+1$, if the feasible setting is the same -- just the first element of the extension tuple is moved to the initial tuple.
We expect the length of the paths determined by the heuristics (and the IP model) to be positively correlated with the value of $c$.

\subsubsection*{Different values of $\Dupl$}\label{subsubsec:ND}
According to the rules of the board game \TT\ $\Dupl = 3$. However, the game is designed for $2$ to $4$ players, we motivates us to consider also smaller numbers of duplicates namely $\Dupl=1,2,3$ in our computational tests. 
A value of $\Dupl$ greater than $1$ implies the possibility that the initial set contains more than one copy of the same node, reducing the number of alternatives for the initialisation. Moreover, the effective value of $c$ reduced also for each iteration of the path extension, since the element of the tuple with be blocked by a duplicate entry. Consequently we assume that the length of the paths is negatively correlated with the value of $\Dupl$.

\subsection{Experiments on the \TT\ graph}\label{subsec:gamegraph}
The graph \G\ has $22$ nodes and $45$ edges (see contains Figure~\ref{Fig:Boardgame as graph}). It is a planar, simple and connected graph which shares some structural similarity with grid graphs.

\subsubsection*{Combination of starting and extension heuristic}\label{subsubsec:normal_allpair}
To evaluate the effectiveness of combinations of initialisation and extension heuristic we test all those combinations on test instances with $\Dupl=2$ and $c=6$, corresponding to the parameter values in the board game \TT. We determine on 10,000 randomly generated $(c,\Dupl)$-settings the resulting paths of all combinations of initialisation and extension heuristic.
The two Tables~\ref{fig:allpair_length} and~\ref{fig:allpair_comptime}) contain the results of this tests. The average node length ranges from $4.3247$ (Start:Random combined with Extension:Random) to $6.1761$ (Start:LongestPath combined with Extension:LongestPath). Thus, in this feasible $(c,\Dupl)$-setting, the combination of the most complex initialisation and extension heuristic yields nearly $50\%$ longer paths on average. The computational time range from $0.002$ seconds (Start:Random combined with Extension:Random) to $0.003$ seconds (Start:LongestPath combined with Extension:LongestPath). Consequently, the combination of the most complex heuristics require roughly $50\%$ more computational time to find a solution as compared to the combination of random heuristics. 
Table~\ref{fig:allpair_length} shows strictly growing average node lengths if the column- or row-index is increased. The only exception of this is between the first and second column of the first row, but their node length difference with a value of $0.0047$ seems to be not significant.

\begin{table}
    \centering
    \small
    \begin{tabular}{l|ccccc}
    \toprule
    Extension: & Random & Degree & Tentacle & Connected & LongestPath\\
     \midrule
    Start:Random  & $4.3247$ & $4.3200$ & $4.3933$ & $4.6493$ & $4.7811$ \\
    Start:Degree  & $4.9768$ & $5.0376$ & $5.1267$ & $5.3814$ & $5.5706$ \\
    Start:Connected  & $5.4374$ & $5.4235$ & $5.5184$ & $5.8933$ & $6.0807$ \\
    Start:LongestPath  & $5.5546$ & $5.5596$ & $5.6294$ & $5.9916$ & $6.1761$ \\\bottomrule
\end{tabular}
    \caption{Average node length of the computed solutions of the comparison of every combination of starting and extension heuristics.}
    \label{fig:allpair_length}
\end{table}

\begin{table}
    \centering
    \small
    \begin{tabular}{l|lllll}
    \toprule
    Extension: & Random & Degree & Tentacle & Connected & LongestPath\\
     \midrule
    Start:Random  & $0.002148$ & $0.002214$ & $0.002469$ & $0.002435$ & $0.002661$ \\
    Start:Degree  & $0.001945$ & $0.002079$ & $0.002347$ & $0.002291$ & $0.002530$ \\
    Start:Connected  & $0.001944$ & $0.002066$ & $0.002338$ & $0.002320$ & $0.002535$ \\
    Start:LongestPath  & $0.002738$ & $0.002836$ & $0.003107$ & $0.003043$ & $0.003314$\\\bottomrule
\end{tabular}
    \caption{Average time (in seconds) used to compute the solutions of the comparison of every combination of starting and extension heuristics.}
    \label{fig:allpair_comptime}
\end{table}



\subsubsection*{Combined heuristics vs.\ IP model} \label{subsubsec:normal_HeuVsIP}
\begin{figure}[htb]
    \centering
    \begin{subfigure}{0.49\textwidth}
    \includegraphics[width=\textwidth]{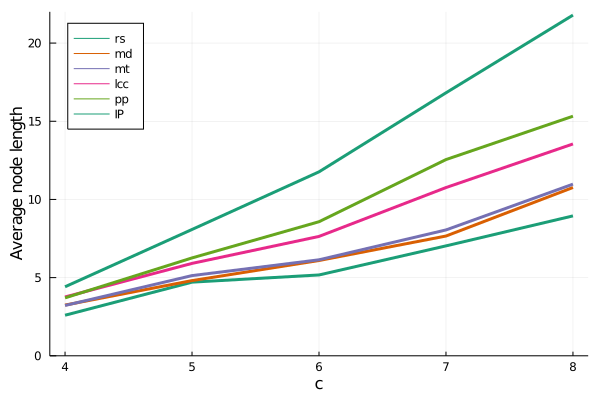}
    \caption{$\Dupl=1$.}
    \end{subfigure}
    \begin{subfigure}{0.49\textwidth}
    \includegraphics[width=\textwidth]{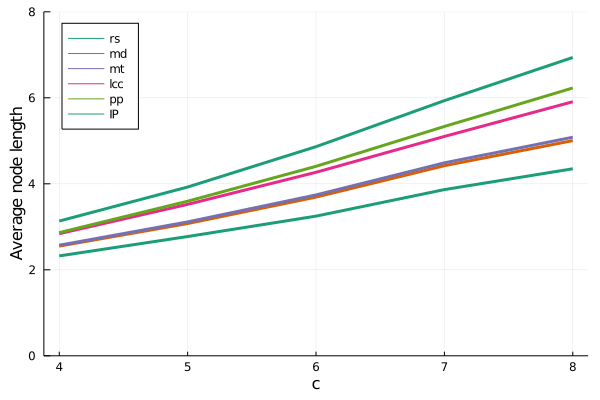}
    \caption{$\Dupl=3$.}
    \end{subfigure}
    \caption{Average node length of the computed solutions on \G. The number of available nodes $c$ is shown on the $x$-axis. \label{Fig:AverageLength_small}}
\end{figure}

While every combinations of an initialisation and an extension heuristic could be applied to (heursitically) solve the online path extension problem, we
will focus in the following on the combinations given in Table~\ref{tab:com_heu}. 
\begin{table}
    \centering
    \small
\begin{tabular}{l|ccccc}
    \toprule
    Extension: & Random & Degree & Tentacle & Connected & LongestPath\\
     \midrule
    Start:Random  & rs &  &  &  &  \\
    Start:Degree  &  & md & mt &  &  \\
    Start:Connected  &  &  &  & lcc &  \\
    Start:LongestPath  &  &  &  &  & pp \\\bottomrule
\end{tabular}
\caption{Considered combinations of initialization and extension heuristics.\label{tab:com_heu}}
\end{table}
The acronyms in Table~\ref{tab:com_heu} have the following meanings: rs -- random search, md -- max degree, mt -- max tentacle, lcc -- largest connected component, pp -- potential path. 
We tied starting and extension heuristics together if they are based on the same idea as introduced in Section~\ref{subsec:heu}. Figure~\ref{Fig:AverageLength_small} shows the average node length for different values of $\Dupl$ and $c$. The average node length of the computed solutions increases for every combined heuristic and the optimal solution of the IP model for increasing values of $c$. However, for increasing values of $\Dupl$ the average determined node length decreases, since the effective value of $c$ is reduced by one each time one node is revealed for the second (or third) time, respectively. If $c$ gets (effectively) smaller, there are less options to choose from. 
The results are based on 100,000 randomly generated $(c,\Dupl)$-settings (for every combination $c=5,\dots ,8$ and $\Dupl=1,2,3$) on which we tested every combined heuristic. The IP model is thereby used to compute the optimal solution based on the complete knowledge. However, the IP model is only solved for $200$ of the $(c,\Dupl)$-settings, due to the high computational effort. 
The average node length of \emph{md} and \emph{mt} are nearly the same for all tested instances on \G . \emph{rs} is the weakest solver since its computed solution yield the smallest average node length. In particular for $\Dupl=2$ and $\Dupl=3$ the solution of the \emph{pp} heuristic reach almost the solution quality of the IP.
The average computational time of all used heuristics is really small. The computational time of the \emph{pp} heuristic increases exponential with the value of $c$ since it solves (multiple) longest path problem (on small instances) in every iteration, but it is at most (in average) $0.009$ seconds per instance. The computational time of all other heuristics only increases a little bit, if $c$ is increased. We assume that the main reason for this increase in computational time is the increasing number of iterations.

\begin{figure}[htb]
    \centering
    \begin{subfigure}{0.75\textwidth}
    \includegraphics[width=\textwidth]{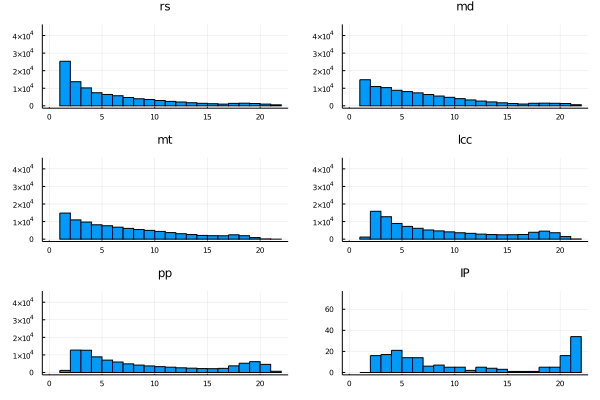}
    \caption{$\Dupl=1$, $c=6$.}
    \end{subfigure}
    
    \vspace{0.5cm}
    \begin{subfigure}{0.75\textwidth}
    \includegraphics[width=\textwidth]{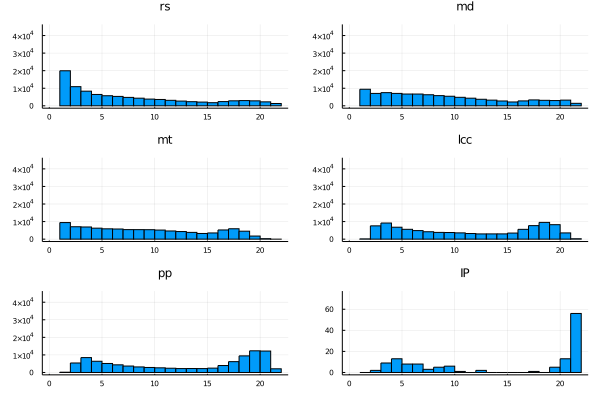}
    \caption{$\Dupl=1$, $c=7$.}
    \end{subfigure}
    \caption{Histograms of the results computed on \G\ for the parameters $\Dupl=1$ and $c=6,7$.\label{Fig:Hist_1_gameboard}}
\end{figure}

\begin{figure}[htb]
    \centering
    \begin{subfigure}{0.75\textwidth}
    \includegraphics[width=\textwidth]{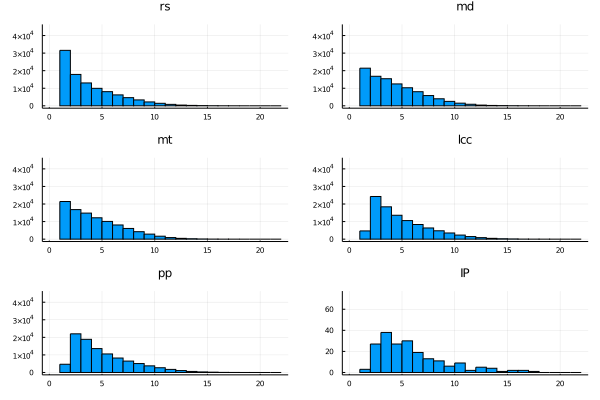}
    \caption{$\Dupl=2$, $c=6$.}
    \end{subfigure}
    
    \vspace{0.5cm}
    \begin{subfigure}{0.75\textwidth}
    \includegraphics[width=\textwidth]{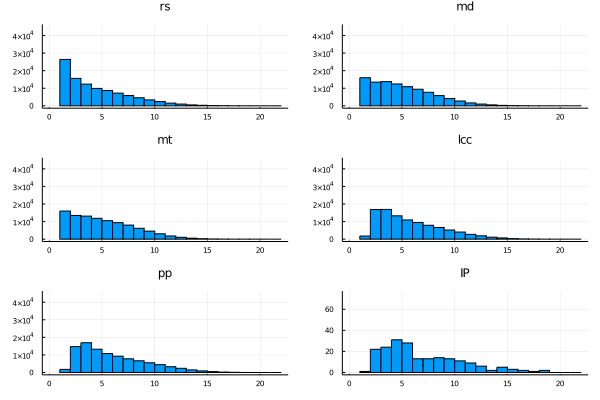}
    \caption{$\Dupl=2$, $c=7$.}
    \end{subfigure}
    \caption{Histograms of the results computed on \G\ for the parameters $\Dupl=2$ and $c=6,7$.\label{Fig:Hist_2_gameboard}}
\end{figure}

\begin{figure}[htb]
    \centering
    \begin{subfigure}{0.75\textwidth}
    \includegraphics[width=\textwidth]{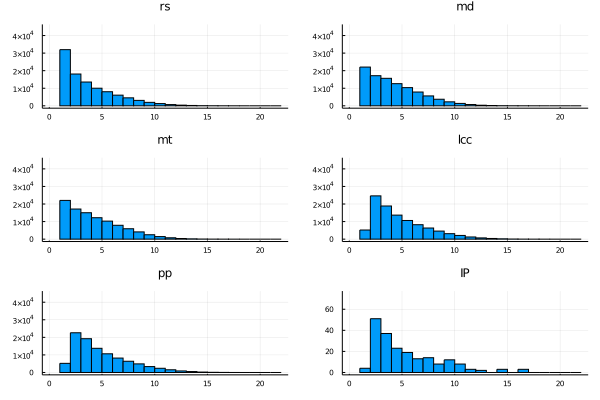}
    \caption{$\Dupl=3$, $c=6$.}
    \end{subfigure}
    
    \vspace{0.5cm}
    \begin{subfigure}{0.75\textwidth}
    \includegraphics[width=\textwidth]{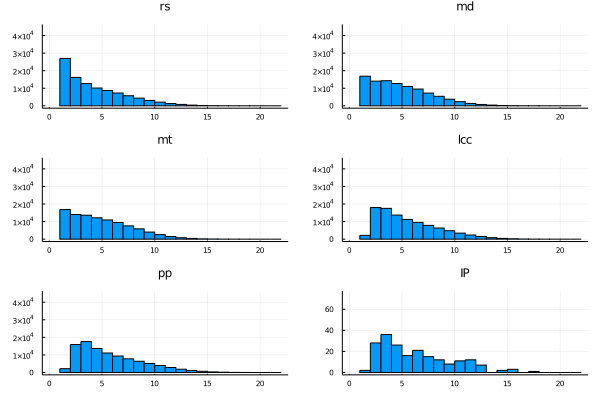}
    \caption{$\Dupl=3$, $c=7$.}
    \end{subfigure}
    \caption{Histograms of the results computed on \G\ for the parameters $\Dupl=3$ and $c=6, 7$.\label{Fig:Hist_3_gameboard}}
\end{figure}

The histograms depicted in Figures~\ref{Fig:Hist_1_gameboard},~\ref{Fig:Hist_2_gameboard} and~\ref{Fig:Hist_3_gameboard} exemplify the detailed results of our computational tests for the parameters $\Dupl=1,2,3$ and $c=6,7$. 
In Figure~\ref{Fig:Hist_1_gameboard} (c) nearly $60$ of the $200$ solutions computed by the IP problem have a node length of $22$ (which is the maximum for this graph) for the $(7,1)$-setting. Apart from the $(7,1)$-setting (and partly the $(6,1)$-setting) all histograms show a similar solution quality for the \emph{pp}-heuristic and the integer program. The results determined by the \emph{lcc}-heuristic achieve a quality almost comparable to \emph{pp}, while it computational time is significantly smaller than \emph{pp}, especially for larger values of $c$. Thus, \emph{lcc} seems to be preferred on larger instances, particularly if the value of $c$ is large.

\begin{figure}[htb]
    \centering
    \begin{subfigure}{0.49\textwidth}
    \includegraphics[width=\textwidth]{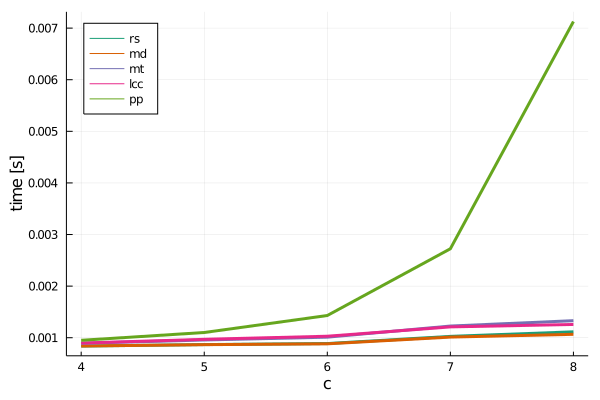}
    \caption{$\Dupl=1$.}
    \end{subfigure}
    \begin{subfigure}{0.49\textwidth}
    \includegraphics[width=\textwidth]{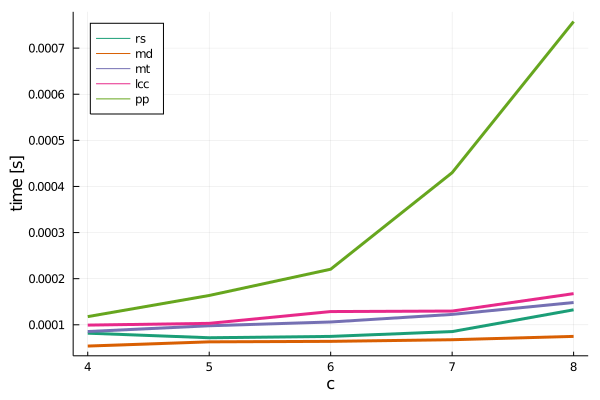}
    \caption{$\Dupl=3$.}
    \end{subfigure}
    \caption{Computational time used by the combined heuristics to compute the solutions on \G. The number of available nodes $c$ is shown on the $x$-axis.\label{Fig:AverageTime_small}}
\end{figure}

\subsubsection*{Tournament of combined heuristics} \label{subsubsec:normal_tournament}

\begin{table}
    \centering
    \small
\begin{tabular}{l|ccccc}
    \toprule
    vs. & rs & md & mt & lcc & pp\\
     \midrule
    rs & $42.0\%$ : $39.9\%$ & $30.5\%$ : $52.8\%$ & $30.4\%$ : $52.9\%$ & $33.0\%$ : $54.2\%$ & $31.1\%$ : $56.1\%$ \\
    md & $56.4\%$ : $27.6\%$ & $44.4\%$ : $35.7\%$ & $44.3\%$ : $35.9\%$ & $46.7\%$ : $38.6\%$ & $44.9\%$ : $40.4\%$ \\
    mt & $56.6\%$ : $27.6\%$ & $44.3\%$ : $35.9\%$ & $44.2\%$ : $35.9\%$ & $46.9\%$ : $38.6\%$ & $44.9\%$ : $40.4\%$ \\
    lcc & $58.7\%$ : $24.8\%$ & $42.1\%$ : $37.9\%$ & $42.1\%$ : $37.9\%$ & $47.1\%$ : $33.9\%$ & $44.3\%$ : $36.4\%$ \\
    pp & $61.1\%$ : $22.9\%$ & $43.7\%$ : $34.6\%$ & $43.7\%$ : $34.6\%$ & $50.4\%$ : $30.6\%$ & $48.4\%$ : $32.5\%$ \\\bottomrule
\end{tabular}
    \caption{Tournament of combined heuristics rs, md, mt, lcc and pp (see Table~\ref{tab:com_heu}). In each entry of the table two relative win frequencies are given. The first is the relative win frequency of the heuristic starting the game (in the same row), the second number is the win frequency of the heuristic having the second move (in the column).\label{tab:tourn}}
\end{table}

We implemented a simple simulation for a two player version of the (OPEP). After the first move of the first player the first node of the extension tuple becomes available and the chosen node is removed from the tuple of available nodes. Then the second player can choose a node from the available nodes to start his own path. After every extension of a path the next node of the extension tuple will be available but the chosen extension node is removed from the available tuple. So, there are exactly $c$ available nodes for every move of both players at any given time, but, if $\Dupl \geq 2$, there can be more than one copy of a node in the tuple of available nodes, lowering the effective number of available nodes the player can choose from. If the value of $\Dupl$ is too small (e.g.~$\Dupl=1$) the extension tuple is too small in a lot of instances. To avoid that, we did our numerical tests with $\Dupl=2$. Each player strictly follows the decisions of one fixed combined heuristic. Note that the used combined heuristics are the same as in the single player tests in the previous section. Every player chooses what is ``best'' in his own understanding of the game, e.g., \emph{md} always chooses the available node with the largest degree and does not adjust the strategy to react to the decisions of the opponent.

If, e.g., the first player is not able to extend the own path at any given iteration, but the second player can extend his own path after the failed attempt of player one, the second player wins this game since his path is (at least) one node longer than the path of the first player. 
If the first player is not able to extend his own path as well as the second player afterwards, the game ends in a tie since both paths have the same node length. We chose $c=6$ to mimic the problem parameters of the board game.

The results of this simulations is shown in table~\ref{tab:tourn}. We tested every combinations of heuristical players on the same set of 100,000 randomly generated $(c,\Dupl)$-settings. The entry of row \emph{md} and colomn \emph{rs} shows the result of the comparison between the combined heuristic \emph{md} (as starting player) and \emph{rs} (as second player): $56.4\%$ : $27.6\%$. This means, that the first player won in $56.4\%$ of the $(c,\Dupl)$-settings and the second player in just $27.6\%$ of settings and consequently $16\%$ of these games ended in a tie.

The starting player always has an advantage if identical combined heuristics are competing against each other (see the diagonal entries of table~\ref{tab:tourn}). If the first player uses \emph{pp}, the best choices for the second player are \emph{md} or \emph{mt}, since the last row of table~\ref{tab:tourn} shows that their win frequency is the highest against a starting \emph{pp}-heuristic. Another observation we want to highlight is the comparison of the two game setting where \emph{rs} and \emph{pp} compete against each other. The win rate of \emph{pp} as starting player is $61.1\%$ compared to $56.1\%$ as the second player. Both values are the highest values for starting and second player among our results, respectively. If \emph{pp} is not the starting player of this setup, its win rate is decreasing by $5\%$, but the win rate of \emph{rs} is increasing by roughly $8\%$, which shows the power of the starting position.

\subsection{Computational tests on a larger graph \Gext}\label{sec:enlargedgraph}
\begin{figure}[htb]
    \small
    \centering
    \begin{tikzpicture}[scale=0.6,every node/.style={draw=black,circle,font=\footnotesize,inner sep=1pt,minimum size=15pt}]
    \node (1) at (0,0) {1};
    \node (2) at (2,0) {2};
    \node (3) at (4,0) {3};
    \node (4) at (6,0) {4};
    \node (5) at (8,0) {5};
    \node (6) at (0,-2) {6};
    \node (7) at (2,-2) {7};
    \node (8) at (4,-2) {8};
    \node (9) at (6,-2) {9};
    \node (10) at (8,-2) {10};
    \node (11) at (-2,-4) {11};
    \node (12) at (0,-4) {12};
    \node (13) at (2,-4) {13};
    \node (14) at (4,-4) {14};
    \node (15) at (6,-4) {15};
    \node (16) at (8,-4) {16};
    \node (17) at (10,-4) {17};
    \node (18) at (-2,-6) {18};
    \node (19) at (0,-6) {19};
    \node (20) at (2,-6) {20};
    \node (21) at (4,-6) {21};
    \node (22) at (6,-6) {22};
    \node[dashed] (23) at (10,0) {23};
    \node[dashed] (24) at (10,-2) {24};
    \node[dashed] (25) at (12,-4) {25};
    \node[dashed] (26) at (-2,-8) {26};
    \node[dashed] (27) at (0,-8) {27};
    \node[dashed] (28) at (2,-8) {28};
    \node[dashed] (29) at (4,-8) {29};
    \node[dashed] (30) at (6,-8) {30};
    \node[dashed] (31) at (8,-6) {31};
    \node[dashed] (32) at (10,-6) {32};
    \draw[-] (1) to (2);
    \draw[-] (1) to (6);
    \draw[-] (1) to (7);
    \draw[-] (2) to (3);
    \draw[-] (2) to (7);
    \draw[-] (3) to (4);
    \draw[-] (3) to (7);
    \draw[-] (3) to (8);
    \draw[-] (3) to (9);
    \draw[-] (4) to (5);
    \draw[-] (4) to (9);
    \draw[-] (4) to (10);
    \draw[-] (6) to (7);
    \draw[-] (6) to (11);
    \draw[-] (7) to (8);
    \draw[-] (7) to (12);
    \draw[-] (7) to (13);
    \draw[-] (8) to (9);
    \draw[-] (8) to (13);
    \draw[-] (8) to (14);
    \draw[-] (8) to (15);
    \draw[-] (9) to (15);
    \draw[-] (9) to (16);
    \draw[-] (10) to (17);
    \draw[-] (11) to (12);
    \draw[-] (11) to (18);
    \draw[-] (11) to (19);
    \draw[-] (12) to (13);
    \draw[-] (12) to (19);
    \draw[-] (12) to (20);
    \draw[-] (13) to (14);
    \draw[-] (13) to (20);
    \draw[-] (14) to (15);
    \draw[-] (14) to (20);
    \draw[-] (14) to (21);
    \draw[-] (15) to (16);
    \draw[-] (15) to (21);
    \draw[-] (15) to (22);
    \draw[-] (16) to (17);
    \draw[-] (16) to (22);
    \draw[-] (17) to (22);
    \draw[-] (18) to (19);
    \draw[-] (19) to (20);
    \draw[-] (20) to (21);
    \draw[-] (21) to (22);
    \draw[dashed] (5) to (23);
    \draw[dashed] (5) to (24);
    \draw[dashed] (17) to (25);
    \draw[dashed] (18) to (26);
    \draw[dashed] (18) to (27);
    \draw[dashed] (19) to (27);
    \draw[dashed] (19) to (28);
    \draw[dashed] (20) to (28);
    \draw[dashed] (21) to (28);
    \draw[dashed] (21) to (29);
    \draw[dashed] (22) to (29);
    \draw[dashed] (22) to (30);
    \draw[dashed] (22) to (31);
    \draw[dashed] (24) to (25);
    \draw[dashed] (25) to (32);
    \draw[dashed] (26) to (27);
    \draw[dashed] (27) to (28);
    \draw[dashed] (28) to (29);
    \draw[dashed] (29) to (30);
    \draw[dashed] (30) to (32);
    \draw[dashed] (31) to (32);
    \end{tikzpicture}
    \caption{Extended version of the board game graph \Gext. Dashed vertices and edges are added to the normal board game graph. This larger graph has $32$ nodes and $66$ edges.\label{fig:largergraph} }
\end{figure}
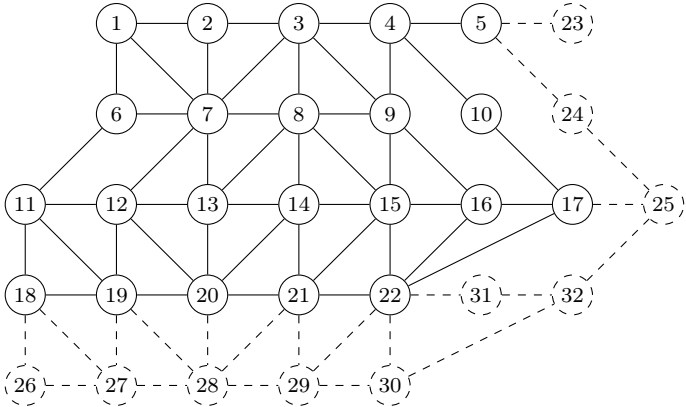

The second part of our computational results was computed for the graph \Gext\ (see Figure~\ref{fig:largergraph}), which contains roughly $50\%$ more nodes and edges than the graph \G\ of \TT, namely $32$ nodes and $66$ edges. We test all combined heuristics and the IP.

\subsubsection*{Results of the combined heuristics}\label{subsubsec:enlarged_HeuVsIP}
We randomly generated 100,000 $(c,\Dupl)$-settings (for every combination of $c=7,\dots ,11$ and $\Dupl=1,2,3$) on which we tested every combined heuristic. The setup of tests is equivalent to the tests done on the smaller graph \G\ with the exception that the value of $c$ is choose $50\%$ larger corresponding to the larger number of node in the graph \Gext.

\begin{figure}[htb]
    \centering
    \begin{subfigure}{0.49\textwidth}
    \includegraphics[width=\textwidth]{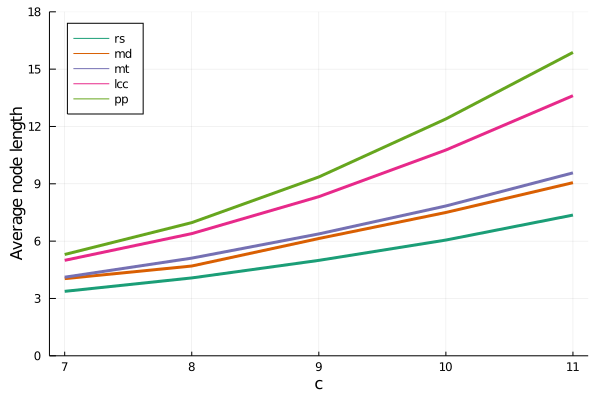}
    \caption{$\Dupl=1$.}
    \end{subfigure}
    \begin{subfigure}{0.49\textwidth}
    \includegraphics[width=\textwidth]{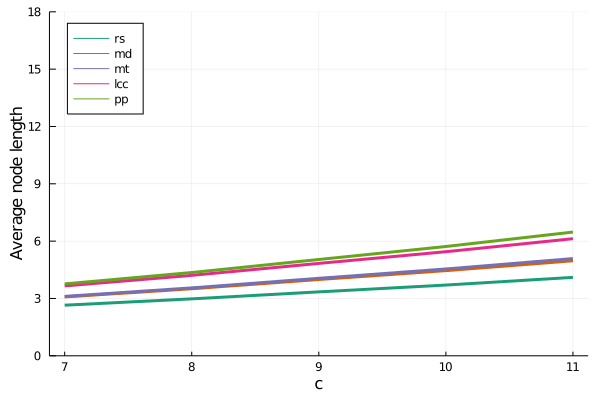}
    \caption{$\Dupl=3$.}
    \end{subfigure}
    \caption{Average length of the computed solutions on the larger artificial graph. The number of available nodes $c$ is shown on the $x$-axis. \label{Fig:length_biggraph}}
\end{figure}

Figure~\ref{Fig:length_biggraph} shows the length of the computed solutions for $\Dupl=1,3$ depending on $c=7,\dots,11$. \emph{pp} yields, in average, the best solutions. The difference in solution quality between different combined heuristics is obvious in Figure~\ref{Fig:length_biggraph} (a), since the paths tend to be longer for smaller values of $\Dupl$. For example, \emph{pp} computes paths that are, in average, double the length of paths computed by \emph{rs}. The overall behavior of the solution quality of all heuristics is similar to the behavior on \G, see~\ref{subsubsec:normal_HeuVsIP}.

The relation of computational times on \G\ is similar to the one \G.  The \emph{pp}-heuristic the computationally most demanding and takes significantly more time for $\Dupl=1$ and larger $c$. The reason, as in the small graph, is the complexity of finding the longest path if the set of vertices grows.

\begin{figure}[htb]
    \centering
    \begin{subfigure}{0.75\textwidth}
    \includegraphics[width=\textwidth]{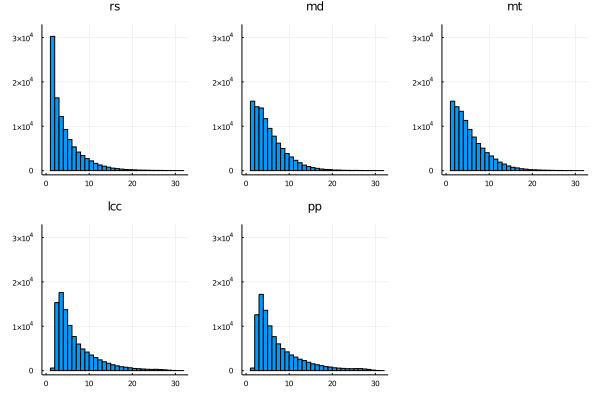}
    \caption{$\Dupl=1$, $c=8$.}
    \end{subfigure}
    
    \vspace*{0.5cm}
    \begin{subfigure}{0.75\textwidth}
    \includegraphics[width=\textwidth]{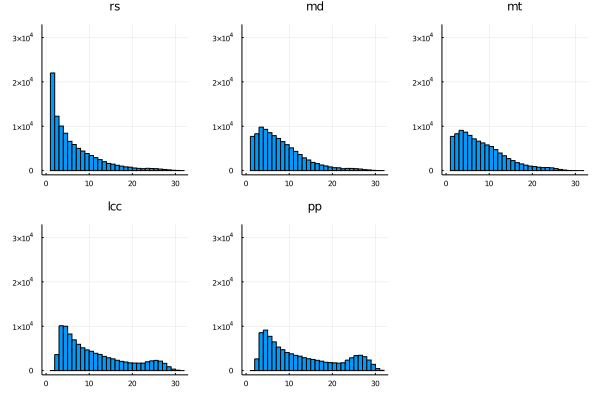}
    \caption{$\Dupl=1$, $c=10$.}
    \end{subfigure}
    \caption{Histograms of the results computed on the larger graph for the parameters $\Dupl=1$ and $c=8, 10$.\label{Fig:Hist_1_large}}
\end{figure}

\begin{figure}[htb]
    \centering
    \begin{subfigure}{0.75\textwidth}
    \includegraphics[width=\textwidth]{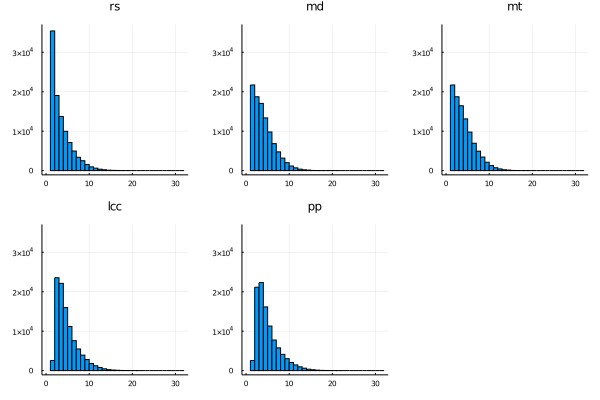}
    \caption{$\Dupl=2$, $c=8$.}
    \end{subfigure}
    
    \vspace*{0.5cm}
    \begin{subfigure}{0.75\textwidth}
    \includegraphics[width=\textwidth]{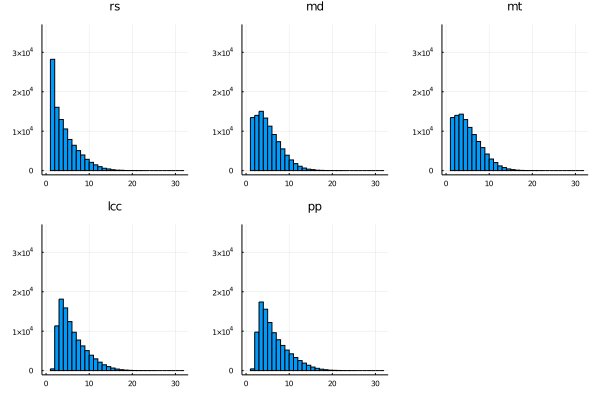}
    \caption{$\Dupl=2$, $c=10$.}
    \end{subfigure}
    \caption{Histograms of the results computed on the larger graph for the parameters $\Dupl=2$ and $c=8, 10$.\label{Fig:Hist_2_large}}
\end{figure}

\begin{figure}[htb]
    \centering
    \begin{subfigure}{0.75\textwidth}
    \includegraphics[width=\textwidth]{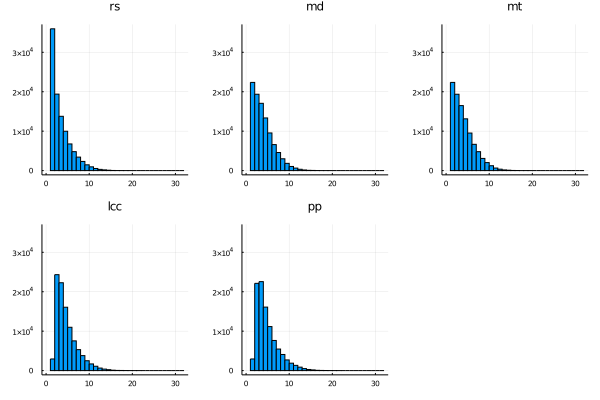}
    \caption{$\Dupl=3$, $c=8$}
    \end{subfigure}
    
    \vspace{0.5cm}
    \begin{subfigure}{0.75\textwidth}
    \includegraphics[width=\textwidth]{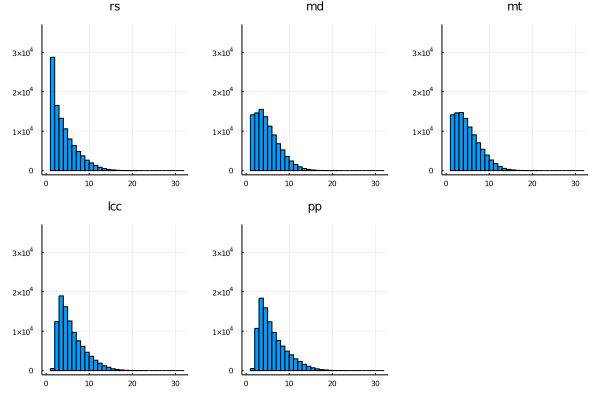}
    \caption{$\Dupl=3$, $c=10$}
    \end{subfigure}
    \\
    \caption{Histograms of the results computed on the larger graph for the parameters $\Dupl=3$ and $c=8, 10$.\label{Fig:Hist_3_large}}
\end{figure}

The different histograms shown in Figure~\ref{Fig:Hist_1_large} (for $\Dupl=1$),~\ref{Fig:Hist_2_large} (for $\Dupl=3$) and~\ref{Fig:Hist_3_large} (for $\Dupl=3$) depict the different quality of solutions using the different combined heuristics. \emph{rs} computes a path of node length one for nearly one third of the $(8,1)$-settings while \emph{pp} obtains nearly no paths shorter with node length less than two. The latter statement holds true for other settings as well. In our numerical tests the heuristics \emph{md} and \emph{mt} obtain on average the same node lengths.

\begin{figure}[htb]
    \centering
    \begin{subfigure}{0.49\textwidth}
    \includegraphics[width=\textwidth]{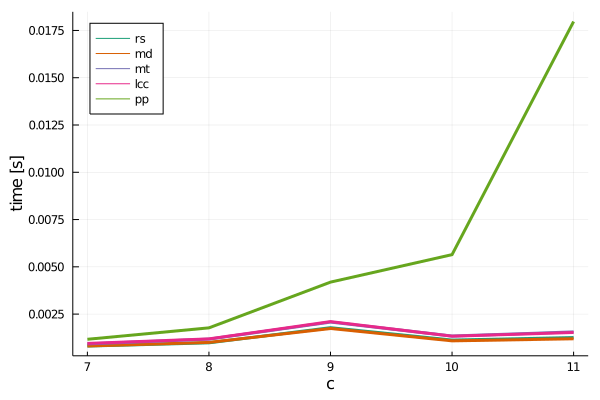}
    \caption{$\Dupl=1$.}
    \end{subfigure}
    \begin{subfigure}{0.49\textwidth}
    \includegraphics[width=\textwidth]{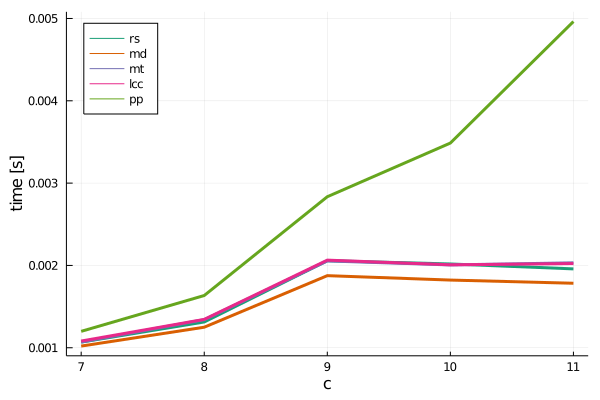}
    \caption{$\Dupl=3$.}
    \end{subfigure}
    \caption{Average Computational time used by the combined heuristics to compute the solutions on the larger artificial graph. The number of available nodes $c$ is shown on the $x$-axis.\label{Fig:time_biggraph}}
\end{figure}

%% file: ConclusionAndOutlook.tex

In this paper, we introduced an online version of a longest path problem that occurs, among others, as a subproblem in the board game \TT. We introduced the concept of tentacles of paths, that is, of nodes that are adjacent to (at least) one of the end nodes of a path while not being nodes in the path. Since paths with many tentacles provide more options for further extensions, extensions to nodes which in turn have many tentacles are preferable. This was confirmed by extensive numerical tests on randomly generated instances on the \TT\ graph.
The numerical results also confirm that heuristics that are based on more involved strategies to predict the potential for further extensions clearly outperform simpler heuristics.

Location and routing problem occur also in many other board games. Prominent examples are, among many others, the games ``The Settlers of Catan'', ``Ticket to Ride'', ``Pandemic'', ``Through the Desert'', ``Mister X'', and many others. Analysing these aspects from a mathematical perspective opens a new research directions that is not only scientifically interesting, but that can also be used for educational purposes to better motivate the underlying mathematical concepts.

%% file: main.bbl
\begin{thebibliography}{10}
\providecommand{\natexlab}[1]{#1}
\providecommand{\url}[1]{\texttt{#1}}
\expandafter\ifx\csname urlstyle\endcsname\relax
  \providecommand{\doi}[1]{doi: #1}\else
  \providecommand{\doi}{doi: \begingroup \urlstyle{rm}\Url}\fi

\bibitem[Ahuja et~al.(1993)Ahuja, Magnanti, and Orlin]{ahuja93network}
Ravindra~K. Ahuja, Thomas~L. Magnanti, and James~B. Orlin.
\newblock \emph{Network Flows: Theory, Algorithms, and Applications}.
\newblock Prentice-Hall, 1993.

\bibitem[Archer et~al.(2011)Archer, Bateni, Hajiaghayi, and Karloff]{Archer}
Aaron Archer, Mohammad~Hossein Bateni, Mohammad~Taghi Hajiaghayi, and Howard
  Karloff.
\newblock Improved approximation algorithms for prize-collecting steiner tree
  and tsp.
\newblock \emph{SIAM Journal on Computing}, 40\penalty0 (2):\penalty0 309--332,
  2011.
\newblock \doi{10.1137/090771429}.

\bibitem[Awerbuch and Kleinberg(2008)]{Awerbuch}
Baruch Awerbuch and Robert Kleinberg.
\newblock Online linear optimization and adaptive routing.
\newblock \emph{Journal of Computer and System Sciences}, 74\penalty0
  (1):\penalty0 97--114, 2008.
\newblock \doi{10.1016/j.jcss.2007.04.016}.

\bibitem[Fernau et~al.(2011)Fernau, Kneis, Kratsch, Langer, Liedloff, Raible,
  and Rossmanith]{Fernau}
Henning Fernau, Joachim Kneis, Dieter Kratsch, Alexander Langer, Mathieu
  Liedloff, Daniel Raible, and Peter Rossmanith.
\newblock An exact algorithm for the maximum leaf spanning tree problem.
\newblock \emph{Theoretical Computer Science}, 412:\penalty0 6290--6302, 2011.
\newblock \doi{10.1016/j.tcs.2011.07.011}.

\bibitem[Fujie(2004)]{Fujie}
Tetsuya Fujie.
\newblock The maximum-leaf spanning tree problem: Formulations and facets.
\newblock \emph{Networks}, 43(4):\penalty0 212--223, 2004.
\newblock \doi{10.1002/net.20001}.

\bibitem[Khabbaz et~al.(2012)Khabbaz, Bhagat, and Lakshmanan]{Khabbaz}
Mohammad Khabbaz, Smriti Bhagat, and Laks V.~S. Lakshmanan.
\newblock Finding heavy paths in graphs: A rank join approach, 2012.

\bibitem[Kneis et~al.(2008)Kneis, Langer, and Rossmanith]{Kneis}
Joachim Kneis, Alexander Langer, and Peter Rossmanith.
\newblock A new algorithm for finding trees with many leaves.
\newblock \emph{Hong SH., Nagamochi H., Fukunaga T. (eds) Algorithms and
  Computation. ISAAC 2008. Lecture Notes in Computer Science, vol 5369.}, 2008.
\newblock \doi{10.1007/978-3-540-92182-0_26}.

\bibitem[Krumke and Noltemeier(2012)]{krumke12graphentheoretische}
Sven~Oliver Krumke and Hartmut Noltemeier.
\newblock \emph{Graphentheoretische Konzepte und Algorithmen}.
\newblock Springer, 2012.
\newblock ISBN 978-3-8348-1849-2.
\newblock \doi{10.1007/978-3-8348-2264-2_13}.

\bibitem[Lu and Ravi(1992)]{Lu}
Hsueh-I Lu and R.~Ravi.
\newblock The power of local optimization: Approximation algorithms for
  maximum-leaf spanning tree.
\newblock In \emph{Proceedings, Thirtieth Annual Allerton Conference on
  Communication, Control and Computing}, pages 533--542, 1992.

\bibitem[Reis et~al.(2015)Reis, Lee, and Usberti]{Reis}
M.~F. Reis, O.~Lee, and F.~L. Usberti.
\newblock Flow-based formulation for the maximum leaf spanning tree problem.
\newblock \emph{Electronic Notes in Discrete Mathematics}, 50:\penalty0
  205--210, 2015.
\newblock \doi{10.1016/j.endm.2015.07.035}.
\newblock LAGOS'15 – VIII Latin-American Algorithms, Graphs and Optimization
  Symposium.

\end{thebibliography}
